%% file: ms.tex
\documentclass[twoside,12pt,a4paper]{article}

\usepackage[latin1]{inputenc}
\usepackage[UKenglish]{babel}
\usepackage{amsmath,amsthm,amscd,amsfonts,amssymb,graphics,graphicx,float}
\usepackage[mathscr]{eucal}
\usepackage{latexsym,cancel,setspace,stmaryrd,newclude,mfirstuc,titlesec,color,dsfont,varwidth,textcomp}
\usepackage[dvipsnames,svgnames,hyperref]{xcolor}
\usepackage[all,color]{xy}
\usepackage{hyperref}
\usepackage[nopostdot]{glossaries}
\usepackage{indentfirst}		
\usepackage{caption}
\usepackage{arydshln}

\usepackage{mathtools}	
\usepackage{multirow}
\usepackage{calc}
\usepackage{upgreek}	

\input{00Preamble}

\title{\bfseries Prym curves with\\
	a vanishing	theta-null}
\author{\normalsize CARLOS MAESTRO PER\'{E}Z}
\date{}

\begin{document}

\maketitle
\thispagestyle{empty}

\begin{abstract}
	\noindent
	If the theta-null divisor $\Theta_{\rm null}$ is moved to the Prym moduli space through the diagram $\K{S}_{g}^{+}\to\K{M}_{g}\leftarrow\K{R}_{g}$, it splits into two irreducible components $\K{P}_{\rm\! null}^{+}$ and $\K{P}_{\rm\! null}^{-}$. Using test curve techniques, we compute the expression of the rational divisor classes of $\OK{P}_{\rm\! null}^{+}$ and $\OK{P}_{\rm\! null}^{-}$ in terms of the generating classes of $\Pic(\OK{R}_{g})_{\B{Q}}$.
\end{abstract}

\input{20CH2}

\bibliographystyle{alpha}
\bibliography{01Ref}

\end{document}

%% file: 00Preamble.tex
\usepackage[stretch=10,shrink=10]{microtype}

\usepackage{enumitem}
\setlist[enumerate]{topsep=-0.5ex,itemsep=0ex,leftmargin=!,labelindent=12mm,labelsep=3mm}
\setlist[enumerate,2]{topsep=0ex,itemsep=0ex,leftmargin=!,labelindent=5mm,labelsep=3mm,label*=\arabic*.}
\setlist[itemize]{topsep=-0.75ex,itemsep=-0.25ex,leftmargin=!,labelindent=12mm,labelsep=3mm}

\colorlet{hip}{teal!100!}			
\colorlet{check}{red!60!black}		
\colorlet{high}{blue!70!black}		
\colorlet{GLOB}{black!}				

\hypersetup{
pdftitle=Prymnull classes,
pdfauthor=Carlos Maestro P\'{e}rez,
pdfsubject=Algebraic Geometry,
bookmarksnumbered=true,bookmarksopen=true,bookmarksopenlevel=2,pdfpagemode=UseOutlines,
colorlinks=true,allcolors=teal}	

\titlespacing*{\section}{0pt}{*5}{*2.5}		
\titlespacing*{\subsection}{0pt}{*5}{2ex}	
\titlespacing*{\subsubsection}{0pt}{*2.5}{1ex}

\input{00Shortcuts}

\usepackage{fancyhdr}
\numberwithin{equation}{section}

\newtheoremstyle{tobeproven}{10pt}{6pt}{\itshape}{}{\bfseries}{.}{ }{}
\theoremstyle{tobeproven}
\newtheorem{thm}[equation]{Theorem}

\newtheorem{lem}[equation]{Lemma}
\newtheorem{prop}[equation]{Proposition}

\newtheoremstyle{statedhere}{10pt}{2pt}{}{}{\bfseries}{.}{ }{}
\theoremstyle{statedhere}
\newtheorem{defn}[equation]{Definition}
\newtheorem{rem}[equation]{Remark}
\newtheorem{exam}[equation]{Example}
\newtheorem*{poss}{Possibility}
\newtheorem*{posss}{Possibilities}

\newtheoremstyle{proofhere}{2pt}{8pt}{}{}{\slshape}{.}{ }{}
\theoremstyle{proofhere}
\newtheorem*{dem}{Proof}

%% file: 00Shortcuts.tex
\newcommand\restr[2]{{
		#1\raisebox{0pt}[0pt][0pt]		{\raisebox{-0.2ex}{$\vert_{#2}$}}}}

\newcommand{\K}[1]{\mathcal{#1}}
\newcommand{\OC}[1]{\smash{\overline{#1}}}
\newcommand{\OK}[1]{\smash{\overline{\mathcal{#1}}}}
\newcommand{\F}[1]{\mathsf{#1}}
\newcommand{\G}[1]{\mathfrak{#1}}
\newcommand{\B}[1]{\mathbb{#1}}
\newcommand{\N}[1]{\mathrm{#1}}

\newcommand{\WT}[1]{\smash{\widetilde{#1}}}
\newcommand{\UC}[1]{\smash{\underline{#1}}}

\DeclareMathOperator{\Pic}{Pic}

\DeclareMathOperator{\Spec}{Spec}

\newcommand{\lto}{\longrightarrow}
\newcommand{\lmto}{\longmapsto}
\newcommand{\fto}{\rightsquigarrow}
\newcommand{\inj}{\hookrightarrow}

%% file: 20CH2.tex
\section{Introduction and preliminaries}
\label{SectCompactMod}

\input{210Compact}

\section{Prym curves and vanishing theta-nulls}
\label{SectPrymPar}

\input{220Dnull}

\section{The geometry of the Prym-null divisors}
\label{SectPnull}

\input{230ClassComp}

%% file: 210Compact.tex
\input{211Introcomp}

\subsection{Stable Prym and spin curves}
\label{SubSectCompStab}

\input{212Stable}

\subsection{Boundary divisors}
\label{SubSectCompBound}

\input{213Piccompact}

\subsection{Test curves}
\label{SubSectCompTest}

\input{214Testcurv}

%% file: 211Introcomp.tex
When looking at the moduli space $\K{S}_{g}^{+}$ of even spin curves of genus $g$, there is a natural, geometric divisor that quickly takes center stage, namely the divisor of curves with a vanishing theta-null, or \textsl{theta-null divisor}:
\[
\begin{array}{rclcl}
	\Theta_{\rm null}&=&\{(C,\theta)\in\K{S}_g^+\;\slash\;h^0(C,\theta)\geq 2\}&\subset&\K{S}_{g}^{+}
\end{array}
\]
In this setting, the theta-null divisor has been thoroughly studied, and its class has been computed by \cite{FarkasEvenSpin}. The same can be said about its pushforward to $\K{M}_{g}$, which was similarly described at an even earlier date in \cite{TeixidorMgnull}. However, its interaction with the moduli space $\K{R}_{g}$ of Prym pairs of genus $g$, by means of the diagram
\[
\xymatrix@R=12pt@C=1pt{
\Theta_{\rm null}&\subset&\K{S}_{g}^{+}
\ar[rrrd]_{\pi_{+}}
&&&&&&\K{R}_{g}\ar[llld]^{\pi_{\K{R}}}
&\supset&\K{P}_{\rm null}&=&(\pi_{\K{R}})^{*}(\pi_{+})_{*}\:\!(\Theta_{\rm null})
\\
&&&&&\K{M}_{g}
}
\]
remains largely unexplored. In this new setting, the divisor $\K{P}_{\rm null}$ splits into two irreducible components
\[
\begin{array}{rcccl}
	\K{P}_{\rm\! null}^{+}&=&\{(C,\eta)\in\K{R}_g\;\slash\;\exists\,\theta\in S_g^+(C)\textrm{ with }\theta\otimes\eta\in\Theta_{\rm null}(C)\}&\subset&\K{R}_{g}
	\\[2mm]
	\K{P}_{\rm\! null}^{-}&=&\{(C,\eta)\in\K{R}_g\;\slash\;\exists\,\theta\in S_g^-(C)\textrm{ with }\theta\otimes\eta\in\Theta_{\rm null}(C)\}&\subset&\K{R}_{g}
\end{array}
\]
which we call \textsl{even} and \textsl{odd Prym-null divisors}. The main result of our work is the computation of the classes of $\OK{P}_{\rm\! null}^{+}$ and $\OK{P}_{\rm\! null}^{-}$ in $\Pic(\OK{R}_{g})_{\B{Q}}$ with standard test curve techniques, culminating in theorem \ref{classprymnulls}. We obtain:
\[
\begin{array}{rccclc}
	\OK{P}_{\rm\! null}^{+}
	&\equiv&\multicolumn{4}{l}{
		2^{g-3}\,\bigg(
		(2^{g-1}+1)\,\lambda
		\hspace*{5.0pt}-\hspace*{5.0pt}
		\dfrac{1}{4}\,\Big(
		2^{g-2}\,\delta_{0}^{\prime}+
		(2^{g-1}+1)\,\delta_{0}^{\rm ram}
		\Big)}
	\vspace*{1.5mm}\\
	&&&-&
	\displaystyle\sum\,
	\Big(
	(2^{i-1}-1)(2^{g-i}-1)\,\delta_{i}+
	(2^{i}-1)(2^{g-i-1}-1)\,\delta_{g-i}
	\hspace*{3.5pt}+
	\vspace*{1mm}\\
	&&&&\multicolumn{1}{r}{+\hspace*{3.5pt}
		(2^{g-1}-2^{i-1}-2^{g-i-1}+1)\,\delta_{i:g-i}
		\Big)\bigg)}
	\vspace*{5mm}\\
	\OK{P}_{\rm\! null}^{-}
	&\equiv&\multicolumn{4}{l}{
		2^{g-3}\,\bigg(
		2^{g-1}\,\lambda
		\hspace*{5.0pt}-\hspace*{5.0pt}
		\dfrac{1}{4}\,\Big(
		2^{g-1}\,\delta_{0}^{\prime\prime}+
		2^{g-2}\,\delta_{0}^{\prime}+
		(2^{g-1}-1)\,\delta_{0}^{\rm ram}
		\Big)}
	\vspace*{1.5mm}\\
	&&&-&
	\displaystyle\sum\,
	\Big(
	2^{i-1}\,(2^{g-i}-1)\,\delta_{i}+
	(2^{i}-1)\,2^{g-i-1}\,\delta_{g-i}
	\hspace*{3.5pt}+
	\vspace*{1mm}\\
	&&&&\multicolumn{1}{r}{+\hspace*{3.5pt}
		(2^{g-1}-2^{i-1}-2^{g-i-1})\,\delta_{i:g-i}
		\Big)\bigg)}
\end{array}
\]
In order to arrive at these expressions, we require some groundwork, which we establish in the remainder of this section. Subsection \ref{SubSectCompStab} offers a brief reminder on the compactifications of $\K{R}_{g}$ and $\K{S}_{g}$, and their corresponding boundaries are described in subsection \ref{SubSectCompBound}. Subsection \ref{SubSectCompTest} introduces the different collections of test curves in $\OK{R}_{g}$ that play a role in the proof of theorem \ref{classprymnulls}. Once the proper background has been set up, the second section deals with the construction of the Prym-null divisors (subsection \ref{SubSectPrymNull}) and their connection with the potential parity change induced on theta characteristics by a Prym root (subsection \ref{SubSectPrymSmooth} for the smooth case, and \ref{SubSectPrymIrred} for the irreducible nodal one). The final section is devoted to the test curve computation, which is carried out in subsections \ref{SubSectPnullRed}, \ref{SubSectPnullEllip} and \ref{SubSectPnullIrred} for reducible nodal curves, curves with elliptic tails and irreducible nodal curves respectively, and whose conclusion leads to subsection \ref{SubSectPnullClass}, where the desired expansions of the even and odd Prym-null classes are given and applied to additional families of Prym curves.

\noindent\textbf{Acknowledgements.} The author would like to thank Prof. Farkas for bringing this problem to his attention, as well as the Berlin Mathematical School and the Humboldt-Universit\"{a}t zu Berlin for their financial support. Special thanks go to Andrei Bud for his helpful discussions regarding remark \ref{infinitecyclic}.

%% file: 212Stable.tex
We work over $\B{C}$, and take $C$ to be a smooth, integral curve of genus $g$.

\begin{defn}
	A \textsl{Prym root} of $C$ is a nontrivial square root of $\K{O}_{C}$, that is, a line bundle $\eta\ncong\K{O}_{C}$ of degree zero equipped with an isomorphism $\eta^{\otimes 2}\cong\K{O}_{C}$.
	A \textsl{Prym pair} is a pair $(C,\eta)$ such that $\eta$ is a Prym root of $C$. The set of Prym roots of $C$ is denoted by $R_g(C)\inj\Pic^{0}(C)-\{\K{O}_C\}$. The moduli space of Prym pairs of genus $g$ is denoted by $\K{R}_{g}$.
\end{defn}

\begin{defn}
	A \textsl{theta characteristic} of $C$ is a square root of $\omega_{C}$, that is, a line bundle $\theta$ of degree $g-1$ equipped with an isomorphism $\theta^{\otimes 2}\cong\omega_{C}$. A \textsl{spin curve} is a pair $(C,\theta)$ such that $\theta$ is a theta characteristic of $C$. The set of theta characteristics of $C$ is denoted by $S_g(C)\inj\Pic^{g-1}(C)$. The moduli space of spin curves of genus $g$ is denoted by $\K{S}_{g}=\K{S}_{g}^{+}\sqcup\K{S}_{g}^{-}$.
\end{defn}

We now want to consider stable versions of these notions. To that end, recall the following standard defitions.

\begin{defn}
	Let $X$ be a complete, connected, nodal curve. We say that $X$ is \textsl{stable} (resp. \textsl{semistable}) if every smooth rational component of $X$ meets the other components of $X$ in at least $3$ points (resp. at least $2$ points).
\end{defn}

\begin{defn}
	Let $E$ be an irreducible component of a semistable curve $X$. Then $E$ is said to be \textsl{exceptional} if it is smooth, rational, and meets the other components in exactly $2$ points.
\end{defn}

\begin{defn}
	Let $X$ be a semistable curve. We say that $X$ is \textsl{quasistable} if any two distinct exceptional components are disjoint. In turn, the \textsl{stable model} of a quasistable curve $X$ is the stable curve $\N{st}(X)$ obtained by contracting each exceptional component to a point.
\end{defn}

\begin{defn}
	A \textsl{stable Prym curve} is a triplet $(X,\eta,\beta)$ where:
	\begin{enumerate}[label=\rm(\roman*)]
		\item $X$ is a quasistable curve (of genus $g$).
		\item $\eta\in\Pic^{0}(X)$ is a nontrivial line bundle of total degree $0$ on $X$ such that $\restr{\eta}{E}=\K{O}_{E}(1)$ for every exceptional component $E$ of $X$.
		\item $\beta\colon\eta^{\otimes 2}\to\K{O}_X$ is a sheaf homomorphism such that the restriction $\restr{\beta}{A}$ is generically non-zero for every non-exceptional component $A$ of $X$.
	\end{enumerate}
	Similarly, a \textsl{stable even spin curve} is a triplet $(X,\theta,\alpha)$ where:
	\begin{enumerate}[label=\rm(\roman*)]
		\item $X$ is a quasistable curve (of genus $g$).
		\item $\theta\in\Pic^{g-1}(X)$ is a line bundle of total degree $g-1$ on $X$ with $h^{0}(X,\theta)$ even, and $\restr{\theta}{E}=\K{O}_{E}(1)$ for every exceptional component $E$ of $X$.
		\item $\alpha\colon\eta^{\otimes 2}\to\omega_X$ is a sheaf homomorphism such that the restriction $\restr{\alpha}{A}$ is generically non-zero for every non-exceptional component $A$ of $X$.
	\end{enumerate}
	For a \textsl{stable odd spin curve}, simply take $h^{0}(X,\theta)$ odd.
\end{defn}

\begin{defn}
	Let $S$ be a scheme. A \textsl{family of stable Prym curves} over the base $S$, or a \textsl{stable Prym curve} over $S$, is a triplet $(f\colon X\to S,\,\eta,\,\beta)$ such that:
	\begin{enumerate}[label=\rm(\roman*)]
		\item $f\colon X\to S$ is a quasistable (genus $g$) curve over $S$.
		\item $\eta\in\Pic^{0}(X)$ is a line bundle on $X$.
		\item $\beta\colon\eta^{\otimes 2}\to\K{O}_X$ is a sheaf homomorphism.
		\item The restriction of $(f\colon X\to S,\,\eta,\,\beta)$ to any fiber $f^{-1}(s)=X_s$ gives rise to a stable Prym curve $(X_s,\eta_s,\beta_{s})$.
	\end{enumerate}
	An isomorphism $(X\to S,\eta,\beta)\cong(X'\to S,\eta',\beta')$ is a pair $(\varphi,\psi)$ where:
	\begin{enumerate}[label=\rm(\roman*)]
		\item $\varphi\colon X\cong X'$ is an isomorphism over $S$.
		\item $\psi\colon\varphi^{*}(\eta')\cong\eta$ is a sheaf isomorphism such that $\varphi^{*}(\beta')=\beta\circ\psi^{\otimes 2}$.
	\end{enumerate}
	With minimal changes, we could likewise define \textsl{families of stable spin curves}.
\end{defn}

The resulting moduli problems all admit proper moduli spaces, namely $\OK{R}_{g}$, $\OK{S}_{g}^{+}$ and $\OK{S}_{g}^{-}$, which respectively compactify $\K{R}_{g}$, $\K{S}_{g}^{+}$ and $\K{S}_{g}^{-}$. Further details on these can be found in \cite{CornalbaCompact}, for the compactification of $\K{S}_{g}$, and \cite{BCFModPrym}, for the Cornalba-inspired compactification of $\K{R}_{g}$. It is also important to highlight \cite{BeauvilleStackyPrym}, where a different but earlier compactification of $\K{R}_{g}$ was built through the use of admissible double covers of stable curves.

Observe that stabilising preserves genus, and consider the natural maps
\[
\begin{array}{rrclccrcl}
	\pi_{\K{R}}\colon&\OK{R}_{g}&\to&\OK{M}_{g},&\quad&
	(X,\eta,\beta)&\mapsto&\N{st}(X)
	\\[1mm]
	\pi_{\K{S}}\colon&\OK{S}_{g}&\to&\OK{M}_{g},&\quad&
	(X,\theta,\alpha)&\mapsto&\N{st}(X)
\end{array}
\]
where $\N{st}(X)$ is the stable model of $X$. These maps are finite and ramified over the boundary, and extend the finite, unramified covers $\K{R}_{g}\to\K{M}_{g}$, $\K{S}_{g}\to\K{M}_{g}$.

%% file: 213Piccompact.tex
In order to study the boundary of $\OK{R}_{g}$, we can take advantage of the map
\[
\begin{array}{rrclccrcl}
	\pi_{\K{R}}\colon&\OK{R}_{g}&\to&\OK{M}_{g},&\quad&
	(X,\eta,\beta)&\mapsto&\N{st}(X)
\end{array}
\]
which turns the decomposition in irreducible components
\[
\begin{array}{ccc}
	\partial\OK{M}_{g}&=&\Delta_{0}\;\cup\;\Delta_{1}\;\cup\;\ldots\;\cup\;\Delta_{\lfloor g/2\rfloor}
\end{array}
\]
into a building block for the corresponding decomposition of $\partial\OK{R}_{g}$.

First, note that a general point $Y\in\Delta_{i}$ is of the form:
\[
\begin{array}{cclccl}
	(i>0)&&Y=C\cup_{p\sim q}D&\;&\textrm{with}&(C,p)\in\K{M}_{i,\,1},\;(D,q)\in\K{M}_{g-i,\,1}
	\\[1mm]
	(i=0)&&Y=B_{pq}&&\textrm{with}&(B,p,q)\in\K{M}_{g-1,\,2}
\end{array}
\]
where $B_{pq}$ denotes the irreducible $1$-nodal curve obtained from $B$ by gluing the points $p$ and $q$. Let us describe the fibers $\pi_{\K{R}}^{-1}(Y)$ for each $i$.

\begin{exam}[$i>0$]\label{BoundDivRgI}
	Let $(X,\eta,\beta)\in\OK{R}_{g}$ with $\N{st}(X)=Y=C\cup_{p\sim q}D$. The existence of $\beta$ prevents $X$ from having exceptional components, i.e.
	\[
	\begin{array}{cccc}
		X=\N{st}(X)=Y=C\cup_{p\sim q}D,&\quad&\beta\colon\eta^{\otimes 2}\cong\K{O}_{Y}=(\K{O}_{C},\K{O}_{D})
	\end{array}
	\]
	Then $\eta$ is a nontrivial element of $J_{2}(C)\oplus J_{2}(D)$, and we have three irreducible components over $\Delta_{i}$, characterized by their general point $(X,\eta,\beta)$:
	\begin{enumerate}[label=\rm(\roman*)]
		\item[$(\Delta_{i}^{\rm n})$]
		\begin{tabular}[t]{|rl}
			Condition:&
			$\eta=(\eta_{C},\K{O}_{D})$ with $\eta_{C}\in R_{i}(C)$.
			\\[1mm]	Notation:&
			$\Delta_{i}^{\rm n}\subset\OK{R}_{g}$ (for \textsl{\UC{n}ontrivial on $i$}), or traditionally $\Delta_{i}$.
			\\[1mm]	Degree:&
			$\deg(\Delta_{i}^{\rm n}\vert\Delta_{i})=2^{2i}-1$.
		\end{tabular}
		\item[$(\Delta_{i}^{\rm t})$]
		\begin{tabular}[t]{|rl}
					Condition:&
			$\eta=(\K{O}_{C},\eta_{D})$ with $\eta_{D}\in R_{g-i}(D)$.
			\\[1mm]	Notation:&
			$\Delta_{i}^{\rm t}\subset\OK{R}_{g}$ (for \textsl{\UC{t}rivial on $i$}), or traditionally $\Delta_{g-i}$.
			\\[1mm]	Degree:&
			$\deg(\Delta_{i}^{\rm t}\vert\Delta_{i})=2^{2(g-i)}-1$.
		\end{tabular}
		\item[$(\Delta_{i}^{\rm p})$]
		\begin{tabular}[t]{|rl}
					Condition:&
			$\eta=(\eta_{C},\eta_{D})$ with $\eta_{C}\in R_{i}(C)$, $\eta_{D}\in R_{g-i}(D)$.
			\\[1mm] Notation:&
			$\Delta_{i}^{\rm p}\subset\OK{R}_{g}$ (for \textsl{\UC{P}rym}), or traditionally $\Delta_{i:g-i}$.
			\\[1mm] Degree:&
			$\deg(\Delta_{i}^{\rm p}\vert\Delta_{i})=(2^{2i}-1)(2^{2(g-i)}-1)$.
		\end{tabular}
	\end{enumerate}
	The pullback of $\Delta_{i}\subset\OK{M}_{g}$ can be written as
	\[
	\pi_{\K{R}}^{*}(\Delta_{i})=\Delta_{i}^{\rm n}+\Delta_{i}^{\rm t}+\Delta_{i}^{\rm p}
	\]
	and, in terms of divisor classes, we have relations
	\[
	\pi_{\K{R}}^{*}(\delta_{i})=\delta_{i}^{\rm n}+\delta_{i}^{\rm t}+\delta_{i}^{\rm p}
	\]
	for $1\leq i\leq\lfloor g/2\rfloor$ and $\delta_{i}^{\rm x}=\K{O}_{\OK{R}_{g}}(\Delta_{i}^{\rm x})\in\Pic(\OK{R}_{g})$, ${\rm x}\in\{{\rm t,n,p}\}$. Observe that
	\[
	\deg(\Delta_{i}^{\rm n}\vert\Delta_{i})+\deg(\Delta_{i}^{\rm t}\vert\Delta_{i})+\deg(\Delta_{i}^{\rm p}\vert\Delta_{i})=2^{2g}-1=\deg(\pi_{\K{R}})
	\]
	as expected.
\end{exam}

\begin{exam}[$i=0$]\label{BoundDivRgO}
	Let $(X,\eta,\beta)\in\OK{R}_{g}$ with $\N{st}(X)=Y=B_{pq}$. There are two possibilities for $X$, depending on whether it contains or not an exceptional component. If it does not, i.e.
	\[
	\begin{array}{cccc}
		X=\N{st}(X)=Y=B_{pq},&\quad&\beta\colon\eta^{\otimes 2}\cong\K{O}_{Y}
	\end{array}
	\]
	then the normalization $\nu\colon B\to B_{pq}$ induces an exact sequence
	\[
	\begin{array}{cccc}
		0\lto\B{Z}_{2}\lto J_{2}(B_{pq})\overset{\nu^{*}}{\lto}J_{2}(B)\lto 0,&\quad&\eta_{B}=\nu^{*}\eta\in J_{2}(B)
	\end{array}
	\]
	and the potential triviality of $\eta_{B}=\nu^{*}\eta$ determines two irreducible components:
	\begin{enumerate}[label=\rm(\roman*)]
		\item[$(\Delta_{0}^{\rm t})$]
		\begin{tabular}[t]{|rl}
			Condition:&
			$\eta_{B}=\K{O}_{B}$, hence $\eta\in(\nu^{*})^{-1}(\eta_{B})-\{\K{O}_{Y}\}$ unique.
			\\[1mm]	Notation:&
			$\Delta_{0}^{\rm t}\subset\OK{R}_{g}$ (for \textsl{\UC{t}rivial}), or traditionally $\Delta_{0}^{\prime\prime}$.
			\\[1mm]	Degree:&
			$\deg(\Delta_{0}^{\rm t}\vert\Delta_{0})=1$.
		\end{tabular}
		\item[$(\Delta_{0}^{\rm p})$]
		\begin{tabular}[t]{|rl}
			Condition:&
			$\eta_{B}\in R_{g-1}(B)$, hence $\eta\in(\nu^{*})^{-1}(\eta_{B})\cong\B{Z}_{2}$.
			\\[1mm]	Notation:&
			$\Delta_{0}^{\rm p}\subset\OK{R}_{g}$ (for \textsl{\UC{P}rym}), or traditionally $\Delta_{0}^{\prime}$.
			\\[1mm]	Degree:&
			$\deg(\Delta_{0}^{\rm p}\vert\Delta_{0})=2\,(2^{2(g-1)}-1)$.
		\end{tabular}
	\end{enumerate}
	On the other hand, if $X$ has an exceptional component $E$, then we can project it onto $Y$ as a sort of ``exceptional blow-up'', i.e. there is a map
	\[
	X=B\cup_{p\sim 0,\,q\sim\infty}E\lto \N{st}(X)=Y=B_{pq}
	\]
	induced by $\nu\colon B\to B_{pq}$, $E\mapsto z=\nu(p)=\nu(q)$. Then we have
	\[
	\begin{array}{cccc}
		\WT{X}=\OC{X-E}\cong B,&\quad&\beta\colon\eta_{B}^{\otimes 2}\cong\K{O}_{B}(-p-q)
	\end{array}
	\]
	for $\eta_{B}=\restr{\eta}{B}\in\Pic(B)$, and Mayer-Vietoris yields an exact sequence
	\[
	\begin{array}{cccc}
		0\lto\B{C}^{*}\lto\Pic(X)\overset{\xi}{\lto}\Pic(B)\oplus\Pic(E)\lto 0,&\quad&\xi(\eta)=(\eta_{B},\K{O}_{E}(1))
	\end{array}
	\]
	This way, we obtain one last irreducible component:
	\begin{enumerate}[label=\rm(\roman*)]
		\item[$(\Delta_{0}^{\rm b})$]
		\begin{tabular}[t]{|rl}
			Condition:&
			$\eta_{B}\in\sqrt{\K{O}_{B}(-p-q)}$.
			\\[1mm]	Notation:&
			$\Delta_{0}^{\rm b}\subset\OK{R}_{g}$ (for \textsl{\UC{b}lown-up}), or traditionally $\Delta_{0}^{{\rm ram}}$.
			\\[1mm]	Degree:&
			$\deg(\Delta_{0}^{\rm b}\vert\Delta_{0})=2^{2(g-1)}$.
		\end{tabular}
	\end{enumerate}
	Due to the appearance of an exceptional component over the node $z\in B_{pq}$, the divisor $\Delta_{0}^{\rm b}$ is in fact the ramification divisor of $\pi_{\K{R}}\colon\OK{R}_{g}\to\OK{M}_{g}$. The pullback of $\Delta_{0}\subset\OK{M}_{g}$ can accordingly be written as
	\[
	\pi_{\K{R}}^{*}(\Delta_{0})=\Delta_{0}^{\rm t}+\Delta_{0}^{\rm p}+2\,\Delta_{0}^{\rm b}
	\]
	and, in terms of divisor classes, we have the relation
	\[
	\pi_{\K{R}}^{*}(\delta_{0})=\delta_{0}^{\rm t}+\delta_{0}^{\rm p}+2\,\delta_{0}^{\rm b}
	\]
	for $\delta_{0}^{\rm x}=\K{O}_{\OK{R}_{g}}(\Delta_{0}^{\rm x})\in\Pic(\OK{R}_{g})$, ${\rm x}\in\{{\rm t,p,b}\}$. Observe that
	\[
	\deg(\Delta_{0}^{\rm t}\vert\Delta_{0})+\deg(\Delta_{0}^{\rm p}\vert\Delta_{0})+2\,\deg(\Delta_{0}^{\rm b}\vert\Delta_{0})=2^{2g}-1=\deg(\pi_{\K{R}})
	\]
	as expected.
\end{exam}

\begin{rem}
	In example \ref{BoundDivRgO}, note that $\deg(\Delta_{0}^{\rm b}\vert\Delta_{0})$ is finite because, for $\eta_{B}\in\sqrt{\K{O}_{B}(-p-q)}$ fixed, any two line bundles
	\[
	\lambda,\mu\in\xi^{-1}(\eta_{B},\K{O}_{E}(1))\cong\B{C}^{*}
	\]
	even if non-isomorphic as bundles, induce triplets
	\[
	(X,\lambda,\beta_{\lambda})\cong(X,\mu,\beta_{\mu})\in\OK{R}_{g}
	\]
	that are always isomorphic as stable Prym curves; see \cite{BCFModPrym} Lemma 2.
\end{rem}

We can now repeat the process for $\OK{S}_{g}$, or rather its irreducible components $\OK{S}_{g}^{+}$, $\OK{S}_{g}^{-}$. Recall the projection
\[
\begin{array}{rrclccrcl}
	\pi_{\K{S}}\colon&\OK{S}_{g}&\to&\OK{M}_{g},&\quad&
	(X,\theta,\alpha)&\mapsto&\N{st}(X)=Y
\end{array}
\]
whose fibers $\pi_{\K{S}}^{-1}(Y)$ we describe for $Y\in\Delta_{i}$ general, $0\leq i\leq\lfloor g/2\rfloor$.

\begin{exam}[$i>0$]\label{BoundDivSgI}
	Let $(X,\theta,\alpha)\in\OK{S}_{g}$ with $\N{st}(X)=Y=C\cup_{p\sim q}D$. The existence of $\alpha$ forces $X$ to have an exceptional component, i.e. there is a map
	\[
	X=C\cup_{p\sim 0}E\cup_{q\sim\infty}D\lto \N{st}(X)=Y=C\cup_{p\sim q}D
	\]
	induced by $E\mapsto z=[p]=[q]$, and we get
	\[
	\begin{array}{cccc}
		\WT{X}=\OC{X-E}\cong C\sqcup D,&\quad&\alpha\colon(\theta_{C},\theta_{D})^{\otimes 2}\cong\restr{\omega_{X}}{\WT{X}}(-p-q)=(\omega_{C},\omega_{D})
	\end{array}
	\]
	for $(\theta_{C},\theta_{D})=\restr{\theta}{\WT{X}}\in\Pic(C)\oplus\Pic(D)$. Therefore, $\theta$ is determined by a pair
	\[
	\begin{array}{cccc}
		(\theta_{C},\theta_{D})\in S_{i}(C)\oplus S_{g-i}(D),&\quad&\theta=(\theta_{C},\K{O}_{E}(1),\theta_{D})\in\Pic(X)
	\end{array}
	\]
	In particular, notice that the (even, odd) parity of $\theta$ is subject to the (identical, alternating) character of the parities of $\theta_{C}$ and $\theta_{D}$, since we have the relation
	\[
	h^{0}(X,\theta)=h^{0}(C,\theta_{C})+h^{0}(D,\theta_{D})
	\]
	by Mayer-Vietoris. As a result, out of the four irreducible components that are obtained over each $\Delta_{i}$, two lie in $\OK{S}_{g}^{+}$ and two lie in $\OK{S}_{g}^{-}$. The even ones are:
	\begin{enumerate}[label=\rm(\roman*)]
		\item[$(\Delta_{i}^{+})$]
		\begin{tabular}[t]{|rl}
			Condition:&
			$\theta_{C}\in S_{i}^{+}(C)$, $\theta_{D}\in S_{g-i}^{+}(D)$.
			\\[1mm]	Notation:&
			$\Delta_{i}^{+}=\Delta_{g-i}^{+}\subset\OK{S}_{g}^{+}$ (for \textsl{even on $i$}), or traditionally $A_{i}^{+}$.
			\\[1mm]	Degree:&
			$\deg(\Delta_{i}^{+}\vert\Delta_{i})=2^{g-1}(2^{i}+1)(2^{g-i}+1)$.
		\end{tabular}
		\item[$(\Delta_{i}^{-})$]
		\begin{tabular}[t]{|rl}
			Condition:&
			$\theta_{C}\in S_{i}^{-}(C)$, $\theta_{D}\in S_{g-i}^{-}(D)$.
			\\[1mm]	Notation:&
			$\Delta_{i}^{-}=\Delta_{g-i}^{-}\subset\OK{S}_{g}^{+}$ (for \textsl{odd on $i$}), or traditionally $B_{i}^{+}$.
			\\[1mm]	Degree:&
			$\deg(\Delta_{i}^{-}\vert\Delta_{i})=2^{g-1}(2^{i}-1)(2^{g-i}-1)$.
		\end{tabular}
	\end{enumerate}
	Similarly, the odd ones are (abusing notation):
	\begin{enumerate}[label=\rm(\roman*)]
		\item[$(\Delta_{i}^{+})$]
		\begin{tabular}[t]{|rl}
			Condition:&
			$\theta_{C}\in S_{i}^{+}(C)$, $\theta_{D}\in S_{g-i}^{-}(D)$.
			\\[1mm]	Notation:&
			$\Delta_{i}^{+}=\Delta_{g-i}^{-}\subset\OK{S}_{g}^{-}$ (for \textsl{even on $i$}), or traditionally $A_{i}^{-}$.
			\\[1mm]	Degree:&
			$\deg(\Delta_{i}^{+}\vert\Delta_{i})=2^{g-1}(2^{i}+1)(2^{g-i}-1)$.
		\end{tabular}
		\item[$(\Delta_{i}^{-})$]
		\begin{tabular}[t]{|rl}
			Condition:&
			$\theta_{C}\in S_{i}^{-}(C)$, $\theta_{D}\in S_{g-i}^{+}(D)$.
			\\[1mm]	Notation:&
			$\Delta_{i}^{-}=\Delta_{g-i}^{+}\subset\OK{S}_{g}^{-}$ (for \textsl{odd on $i$}), or traditionally $B_{i}^{-}$.
			\\[1mm]	Degree:&
			$\deg(\Delta_{i}^{-}\vert\Delta_{i})=2^{g-1}(2^{i}-1)(2^{g-i}+1)$.
		\end{tabular}
	\end{enumerate}
	Observe that a factor of $2$ has to be considered in the computation
	\[
	\begin{array}{ccccl}
		\Delta_{i}^{+}\subset\OK{S}_{g}^{+},&&
		\deg(\Delta_{i}^{+}\vert\Delta_{i})&=&
		2\cdot\#S_{i}^{+}(C)\cdot\#S_{g-i}^{+}(D)
		\\[1mm]
		&&&=&2\cdot 2^{i-1}(2^{i}+1)\cdot 2^{g-i-1}(2^{g-i}+1)
		\\[1mm]
		&&&=&2\cdot 2^{g-2}(2^{i}+1)(2^{g-i}+1)
	\end{array}
	\]
	to account for the nontrivial automorphism of $(X,\theta,\alpha)$ that arises from scaling by $-1$ on the exceptional component. Over the coarse moduli space $\OC{S}_{g}^{+}\to\OC{M}_{g}$, this factor is not present. Consequently, the pullback of $\Delta_{i}\subset\OK{M}_{g}$ (resp. $\subset\OC{M}_{g})$ by $\pi_{+}\colon\OK{S}_{g}^{+}\to\OK{M}_{g}$ (resp. $\pi_{+}\colon\OC{S}_{g}^{+}\to\OC{M}_{g}$) can be written as
	\[
	\begin{array}{ccc}
		\pi_{+}^{*}(\Delta_{i})=2\,\Delta_{i}^{+}+2\,\Delta_{i}^{-}&\quad&\textrm{(resp. $\pi_{+}^{*}(\Delta_{i})=\Delta_{i}^{+}+\Delta_{i}^{-}$)}
	\end{array}
	\]
	and, in terms of divisor classes, we have relations
	\[
	\begin{array}{ccc}
		\;\;\pi_{+}^{*}(\delta_{i})=2\,\delta_{i}^{+}+2\,\delta_{i}^{-}&\quad&\textrm{(resp. $\pi_{+}^{*}[\Delta_{i}]=[\Delta_{i}^{+}]+[\Delta_{i}^{-}]$)}
	\end{array}
	\]
	for $1\leq i\leq\lfloor g/2\rfloor$ and $\delta_{i}^{\rm x}=\K{O}_{\OK{S}_{g}^{+}}(\Delta_{i}^{\rm x})\in\Pic(\OK{S}_{g}^{+})$, ${\rm x}\in\{+,-\}$. The same analysis works for $\OK{S}_{g}^{-}$ and $\OC{S}_{g}^{-}$.
\end{exam}

\begin{exam}[$i=0$]\label{BoundDivSgO}
	Let $(X,\theta,\alpha)\in\OK{S}_{g}$ with $\N{st}(X)=Y=B_{pq}$. There are again two possibilities for $X$. If it has no exceptional components, i.e.
	\[
	\begin{array}{cccc}
		X=\N{st}(X)=Y=B_{pq},&\quad&\alpha\colon\theta^{\otimes 2}\cong\omega_{Y}
	\end{array}
	\]
	then the normalization $\nu\colon B\to B_{pq}$ induces a double cover
	\[
	\begin{array}{cccc}
		\nu^{*}\colon\sqrt{\omega_{Y}\vphantom{(p+q)}}\lto\sqrt{\omega_{B}(p+q)},&\quad&(\nu^{*}\theta)^{\otimes 2}\cong\nu^{*}\omega_{Y}\cong\omega_{B}(p+q)
	\end{array}
	\]
	so that $\theta$ is determined by a square root $\theta_{B}\in\sqrt{\omega_{B}(p+q)}$ and a choice on how to glue its fibers $\restr{\theta_{B}}{p}$ and $\restr{\theta_{B}}{q}$. Only two such gluings are possible, one making $h^{0}(X,\theta)$ even and the other one making it odd. We describe the component $\Delta_{0}^{\rm n}$ obtained in this way only for $\OK{S}_{g}^{+}$, as its $\OK{S}_{g}^{-}$ counterpart is very similar.
	\begin{enumerate}[label=\rm(\roman*)]
		\item[$(\Delta_{0}^{\rm n})$]
		\begin{tabular}[t]{|rl}
			Condition:&
			$\theta_{B}\in\sqrt{\omega_{B}(p+q)}$ with even gluing.
			\\[1mm]	Notation:&
			$\Delta_{0}^{\rm n}\subset\OK{S}_{g}^{+}$ (for \textsl{\UC{n}ot blown-up}), or traditionally $A_{0}^{+}$.
			\\[1mm]	Degree:&
			$\deg(\Delta_{0}^{\rm n}\vert\Delta_{0})=2^{2g-2}$.
		\end{tabular}
	\end{enumerate}
	On the other hand, if $X$ has an exceptional component $E$, then
	\[
	\begin{array}{cccc}
		X=B\cup_{p\sim 0,\,q\sim\infty}E\lto \N{st}(X)=Y=B_{pq},&\quad&
		\WT{X}=\OC{X-E}\cong B
	\end{array}
	\]
	and we have $\theta_{B}=\restr{\theta}{B}\in S_{g-1}(B)$, since $\alpha$ gives rise to an isomorphism
	\[
	\alpha\colon\theta_{B}^{\otimes 2}\cong\restr{\omega_{X}}{B}(-p-q)\cong\nu^{*}\omega_{Y}(-p-q)\cong\omega_{B}
	\]
	Moreover, recall the exact sequence
	\[
	\begin{array}{cccc}
		0\lto\B{C}^{*}\lto\Pic(X)\overset{\xi}{\lto}\Pic(B)\oplus\Pic(E)\lto 0,&\quad&\xi(\theta)=(\theta_{B},\K{O}_{E}(1))
	\end{array}
	\]
	and note that $h^{0}(X,\theta)=h^{0}(B,\theta_{B})$, again by Mayer-Vietoris. In conclusion, we get the remaining irreducible component of $\partial\OK{S}_{g}^{+}$, and similarly for $\partial\OK{S}_{g}^{-}$, as:
	\begin{enumerate}[label=\rm(\roman*)]
		\item[$(\Delta_{0}^{\rm b})$]
		\begin{tabular}[t]{|rl}
			Condition:&
			$\theta_{B}\in S_{g-1}^{+}(B)$.
			\\[1mm]	Notation:&
			$\Delta_{0}^{\rm b}\subset\OK{S}_{g}^{+}$ (for \textsl{\UC{b}lown-up}), or traditionally $B_{0}^{+}$.
			\\[1mm]	Degree:&
			$\deg(\Delta_{0}^{\rm b}\vert\Delta_{0})=2^{g-2}(2^{g-1}+1)$.
		\end{tabular}
	\end{enumerate}
	The pullback of $\Delta_{0}\subset\OK{M}_{g}$ by $\pi_{+}\colon\OK{S}_{g}^{+}\to\OK{M}_{g}$ can be written as
	\[
	\pi_{+}^{*}(\Delta_{0})=\Delta_{0}^{\rm n}+2\,\Delta_{0}^{\rm b}
	\]
	and, in terms of divisor classes, we have the relation
	\[
	\pi_{+}^{*}(\delta_{0})=\delta_{0}^{\rm n}+2\,\delta_{0}^{\rm b}
	\]
	for $\delta_{0}^{\rm x}=\K{O}_{\OK{S}_{g}^{+}}(\Delta_{0}^{\rm x})\in\Pic(\OK{S}_{g}^{+})$, ${\rm x}\in\{{\rm n,b}\}$. Finally, the divisor
	\[
	\begin{array}{ccc}
	 	\Delta_{0}^{\rm b}+\sum(\Delta_{i}^{+}+\Delta_{i}^{-})&\quad&\textrm{(resp. $\Delta_{0}^{\rm b}$)}
	\end{array}
	\]
	is the ramification divisor of $\pi_{+}\colon\OK{S}_{g}^{+}\to\OK{M}_{g}$ (resp. $\pi_{+}\colon\OC{S}_{g}^{+}\to\OC{M}_{g}$).
\end{exam}

Examples \ref{BoundDivRgI} and \ref{BoundDivRgO} provide us with a collection of boundary classes of $\OK{R}_{g}$, while examples \ref{BoundDivSgI} and \ref{BoundDivSgO} follow suit with $\OK{S}_{g}^{+}$ (and $\OK{S}_{g}^{-}$):
\[
\begin{array}{ccclcc}
	\delta_{0}^{\rm t},\,\delta_{0}^{\rm p},\,\delta_{0}^{\rm b},&\delta_{i}^{\rm t},\,\delta_{i}^{\rm n},\,\delta_{i}^{\rm p}&\in&\Pic(\OK{R}_{g}),&\quad&1\leq i\leq\lfloor g/2\rfloor
	\\[2mm]
	\delta_{0}^{\rm n},\,\delta_{0}^{\rm b},&\delta_{i}^{+},\,\delta_{i}^{-}&\in&\Pic(\OK{S}_{g}^{+}),&\quad&1\leq i\leq\lfloor g/2\rfloor
\end{array}
\]
For $g\geq 5$, we then get
\[
\begin{array}{ccc}
	\Pic(\OK{R}_{g})_{\B{Q}}&=&\displaystyle\lambda\,\B{Q}\;\oplus\;\delta_{0}^{\rm t}\,\B{Q}\;\oplus\;\delta_{0}^{\rm p}\,\B{Q}\;\oplus\;\delta_{0}^{\rm b}\,\B{Q}\;\oplus\;\bigoplus_{i=1}^{\lfloor g/2\rfloor}(\delta_{i}^{\rm t}\,\B{Q}\;\oplus\;\delta_{i}^{\rm n}\,\B{Q}\;\oplus\;\delta_{i}^{\rm p}\,\B{Q})
\end{array}
\]
and similarly
\[
\begin{array}{ccc}
	\Pic(\OK{S}_{g}^{+})_{\B{Q}}&=&\displaystyle\lambda\,\B{Q}\;\oplus\;\delta_{0}^{\rm n}\,\B{Q}\;\oplus\;\delta_{0}^{\rm b}\,\B{Q}\;\oplus\;\bigoplus_{i=1}^{\lfloor g/2\rfloor}(\delta_{i}^{+}\,\B{Q}\;\oplus\;\delta_{i}^{-}\,\B{Q})
\end{array}
\]
where $\lambda$ denotes the pullback of $\lambda\in\Pic(\OK{M}_{g})$ to $\OK{R}_{g}$ and $\OK{S}_{g}^{+}$, respectively.

\begin{rem}
	As the Hodge bundle construction used to build $\lambda\in\Pic(\OK{M}_{g})$ commutes with base change, the class $\lambda$ in $\OK{R}_{g}$ or $\OK{S}_{g}^{+}$ can likewise be defined by means of the Hodge bundle associated to each of these spaces.
\end{rem}

\begin{rem}\label{infinitecyclic}
	For the decompositions of $\Pic(\OK{R}_{g})_{\B{Q}}$ and $\Pic(\OK{S}_{g}^{+})_{\B{Q}}$ to hold, it is enough to see that $\Pic(\K{R}_{g})_{\B{Q}}$ and $\Pic(\K{S}_{g}^{+})_{\B{Q}}$ are infinite cyclic. In the case of $\K{R}_{g}$, we have finite maps
	\[
	\K{M}_{g}(2)\to\K{R}_{g}\to\K{M}_{g},\qquad(C,\eta_{1},\ldots,\eta_{g})\mapsto(C,\eta_{1})\mapsto C
	\]
	where $\K{M}_{g}(2)$ is the moduli of curves with a level $2$ structure, that is, a basis of the $2$-torsion of their Jacobian. As a result, we get injective pullback maps
	\[
	\Pic(\K{M}_{g})_{\B{Q}}\inj\Pic(\K{R}_{g})_{\B{Q}}\inj\Pic(\K{M}_{g}(2))_{\B{Q}}
	\]
	Since Putman's work \cite{PutmanPicLevel, PutmanSecondHom} shows that $\Pic(\K{M}_{g}(2))_{\B{Q}}\cong\B{Q}$ for $g\geq 5$, it follows that $\Pic(\K{R}_{g})_{\B{Q}}\cong\B{Q}$ in this range. The corresponding result for $\K{S}_{g}^{+}$ is due to Harer \cite{HarerPicSpin} for $g\geq 9$, and Putman \cite{PutmanPicLevel} for $g\geq 5$.
\end{rem}

%% file: 214Testcurv.tex
Let us take the most basic families of test curves on $\OK{M}_{g}$ and examine ways of lifting them to $\OK{R}_{g}$. In the following examples, we denote
\[
\pi\colon\OK{R}_{g}\to\OK{M}_{g}
\]
instead of $\pi_{\K{R}}$, as we will not work with spin curves here. However, descriptions of common test curves on $\OK{S}_{g}$ can be found in \cite{FarkasEvenSpin} or \cite{FVOddSpin}.

\begin{exam}[reducible nodal curves]\label{testrednod}
	For each integer $2\leq i\leq g-1$, we fix general curves $C\in\K{M}_i$ and $(D,q)\in\K{M}_{g-i,\,1}$ and consider the test curve
	\[
	\K{C}^{i}=(C\times C)\cup_{\Delta_{C}\:\!\sim\, C\times\{q\}}(C\times D)\lto C
	\]
	corresponding to the family of reducible nodal curves
	\[
	\begin{array}{rclclcl}
		\K{C}^{i}&\equiv&\{C\cup_{y\sim q}D\}_{y\in C}&\subset&\Delta_i&\subset&\OK{M}_g
	\end{array}
	\]
	Using the standard test curve techniques of \cite{HMModuliCurves} Chapter 3, we can see that the intersection numbers of $\K{C}^{i}$ with the generators of $\Pic(\OK{M}_g)_{\B{Q}}$ given earlier in the section are described by the table:
	\[
	\begin{array}{c|cccccccc}
		&\lambda&\delta_{i}&\delta_{(j\neq i)}\\[\dimexpr-\normalbaselineskip+0.5mm]
		\\\hline
		\\[\dimexpr-\normalbaselineskip+1mm]
		\K{C}^i&0&2-2i&0
	\end{array}
	\]
	We now fix two Prym roots $\eta_C\in R_i(C)$, $\eta_D\in R_{g-i}(D)$ and lift $\K{C}^{i}$ to test curves $F_i$, $G_i$, $H_i\to C$, as follows:
	\[
	\begin{array}{rccclcl}
		F_i&\equiv&\{(C\cup_{y\sim q}D,\,(\eta_C,\K{O}_D))\}_{y\in C}&\subset&\Delta_{i}^{\rm n}&\subset&\OK{R}_g
		\\[2mm]
		G_i&\equiv&\{(C\cup_{y\sim q}D,\,(\K{O}_C,\eta_D))\}_{y\in C}&\subset&\Delta_{i}^{\rm t}&\subset&\OK{R}_g
		\\[2mm]
		H_i&\equiv&\{(C\cup_{y\sim q}D,\,(\eta_C,\eta_D))\}_{y\in C}&\subset&\Delta_{i}^{\rm p}&\subset&\OK{R}_g
	\end{array}
	\]
	Observe that $\pi_{*}(F_i)=\pi_{*}(G_i)=\pi_{*}(H_i)=\K{C}^{i}$. Then
	\[\begin{array}{ccccccl}
		F_i\cdot\delta_{i}^{\rm n}&=&F_i\cdot\pi^{*}\delta_{i}&=&\K{C}^{i}\cdot\delta_{i}&=&2-2i
		\\[2mm]
		G_i\cdot\delta_{i}^{\rm t}&=&G_i\cdot\pi^{*}\delta_{i}&=&\K{C}^{i}\cdot\delta_{i}&=&2-2i
		\\[2mm]
		H_i\cdot\delta_{i}^{\rm p}&=&H_i\cdot\pi^{*}\delta_{i}&=&\K{C}^{i}\cdot\delta_{i}&=&2-2i
	\end{array}
	\]
	and all other intersection numbers are $0$, which is collected in the table:
	\[
	\begin{array}{c|cccccccc}
		&\lambda&\delta_{0}^{\rm t}&\delta_{0}^{\rm p}&\delta_{0}^{\rm b}&\delta_{i}^{\rm n}&\delta_{i}^{\rm t}&\delta_{i}^{\rm p}&\delta_{(j\neq i)}\\[\dimexpr-\normalbaselineskip+0.5mm]
		\\\hline
		\\[\dimexpr-\normalbaselineskip+1mm]
		F_i&0&0&0&0&2-2i&0&0&0\\
		G_i&0&0&0&0&0&2-2i&0&0\\
		H_i&0&0&0&0&0&0&2-2i&0
	\end{array}
	\]
	Note the exception of $g=2i$, where we have $F_i\cdot\delta_{i}^{\rm n}=G_i\cdot\delta_{i}^{\rm n}=2-2i$.
\end{exam}

Let us move on to the standard degree $12$ pencil of elliptic tails in $\OK{M}_g$.

\begin{exam}[elliptic tails]\label{testelliptail}
	We fix a general curve $(C,p)\in\K{M}_{g-1,\,1}$ and a general pencil $f\colon\N{Bl}_{9}(\B{P}^2)\to\B{P}^{1}$ of plane cubics, with fibers
	\[
	\begin{array}{rcl}
		\{E_{\lambda}=f^{-1}(\lambda)\}_{\lambda\in\B{P}^1}&\subset&\OK{M}_{1}
	\end{array}
	\]
	together with a section $\sigma\colon\B{P}^{1}\to\N{Bl}_{9}(\B{P}^2)$ induced by one of the basepoints. We may then glue the curve $(C,p)$ to the pencil $f$ along $\sigma$, thus producing a pencil of stable curves
	\[
	\K{C}^{0}=(C\times \B{P}_{1})\cup_{\{p\}\times\B{P}^{1}\:\!\sim\, \sigma(\B{P}^{1})}\N{Bl}_{9}(\B{P}^2)\lto \B{P}^{1}
	\]
	which corresponds to
	\[
	\begin{array}{rclclcl}
		\K{C}^{0}&\equiv&\{C\cup_{p\sim\sigma(\lambda)}E_{\lambda}\}_{\lambda\in\B{P}^{1}}&\subset&\Delta_{1}&\subset&\OK{M}_{g}
	\end{array}
	\]
	As in the previous example, \cite{HMModuliCurves} shows that the intersection numbers of the pencil $\K{C}^{0}$ with the generators of $\Pic(\OK{M}_g)_{\B{Q}}$ are given by the table:
	\[
	\begin{array}{c|cccccccc}
		&\lambda&\delta_{0}&\delta_{1}&\delta_{(j\geq 2)}\\[\dimexpr-\normalbaselineskip+0.5mm]
		\\\hline
		\\[\dimexpr-\normalbaselineskip+1mm]
		\K{C}^0&1&12&-1&0
	\end{array}
	\]
	If we now fix a Prym root $\eta_C\in R_{g-1}(C)$, then the degree $3$ branched covering
	\[
	\gamma_1\colon\OK{R}_{1,1}\to\OK{M}_{1,1}
	\]
	allows us to lift $\K{C}^{0}$ to test curves $F_0$, $G_0$, $H_0$, as follows:
	\[
	\begin{array}{rclclcl}
		F_0&\equiv&\{(C\cup_{p\sim\sigma(\lambda)}E_{\lambda},\,(\eta_C,\K{O}_{E_{\lambda}}))\}_{\lambda\in\B{P}^1}&\subset&\Delta_{1}^{\rm t}&\subset&\OK{R}_g
		\\[2mm]
		G_0&\equiv&\{(C\cup_{p\sim\sigma(\lambda)}E_{\lambda},\,(\K{O}_C,\eta_{E_{\lambda}}))\;\slash\;\eta_{E_{\lambda}}\in\gamma_1^{-1}(E_{\lambda})\}_{\lambda\in\B{P}^1}&\subset&\Delta_{1}^{\rm n}&\subset&\OK{R}_g
		\\[2mm]
		H_0&\equiv&\{(C\cup_{p\sim\sigma(\lambda)}E_{\lambda},\,(\eta_C,\eta_{E_{\lambda}}))\;\slash\;\eta_{E_{\lambda}}\in\gamma_1^{-1}(E_{\lambda})\}_{\lambda\in\B{P}^1}&\subset&\Delta_{1}^{\rm p}&\subset&\OK{R}_g
	\end{array}
	\]
	Observe that $\pi_{*}(F_0)=\K{C}^{0}$ and $\pi_{*}(G_0)=\pi_{*}(H_0)=3\,\K{C}^{0}$, so in particular
	\[
	\begin{array}{ccccccl}
		F_0\cdot\delta_{1}^{\rm t}&=&F_0\cdot\pi^{*}\delta_{1}&=&\K{C}^{0}\cdot\delta_{1}&=&-1
		\\[2mm]
		G_0\cdot\delta_{1}^{\rm n}&=&G_0\cdot\pi^{*}\delta_{1}&=&3\,\K{C}^{0}\cdot\delta_{1}&=&-3
		\\[2mm]
		H_0\cdot\delta_{1}^{\rm p}&=&H_0\cdot\pi^{*}\delta_{1}&=&3\,\K{C}^{0}\cdot\delta_{1}&=&-3
	\end{array}
	\]
	Looking at the $12$ points $\lambda_{\infty}\in\B{P}^1$ that correspond to singular fibers of $\K{C}^{0}$ and blowing up the node of the rational component $E_{\lambda_{\infty}}\in\Delta_{0}$, we see that, for $F_0$, the pullback of $\eta_{\lambda_{\infty}}=(\eta_C,\K{O}_{E_{\lambda_{\infty}}})$ is $(\eta_C,\K{O}_{\B{P}^1})$, which is nontrivial. As discussed in example \ref{BoundDivRgO}, this implies that $F_{0,\,\lambda_{\infty}}\in\Delta_{0}^{\rm p}$, hence
	\[
	\begin{array}{ccccccc}
		F_0\cdot \delta_{0}^{\rm p}&=&F_0\cdot\pi^{*}\delta_{0}&=&\K{C}^{0}\cdot\delta_{0}&=&12
	\end{array}
	\]
	Furthermore, the covering $\gamma_1\colon\OK{R}_{1,1}\to\OK{M}_{1,1}$ is branched over $E_{\lambda_{\infty}}$, and thus the fiber $\gamma_1^{-1}(E_{\lambda_{\infty}})$ consists of two elements: one lying in the ramification divisor of $\gamma_1$, which we denote by $\eta_{E_{\lambda_{\infty}}}^{\rm b}$, and one outside, which we denote by $\eta_{E_{\lambda_{\infty}}}^{\rm t}$. Then the pullback of $(\K{O}_C,\eta_{E_{\lambda_{\infty}}}^{\rm t})$ is $(\K{O}_C,\K{O}_{\B{P}^1})$, that is, $(\K{O}_C,\eta_{E_{\lambda_{\infty}}}^{\rm t})\in\Delta_{0}^{\rm t}$, and we get
	\[
	\begin{array}{ccccl}
		G_0\cdot\delta_{0}^{\rm t}&=&\K{C}^{0}\cdot\delta_{0}&=&12
		\\[2mm]
		G_0\cdot\delta_{0}^{\rm b}&=&\K{C}^{0}\cdot\delta_{0}&=&12
	\end{array}
	\]
	Finally, the pair $(\eta_C,\eta_{E_{\lambda_{\infty}}}^{\rm t})$ pulls back to the nontrivial pair $(\eta_C,\K{O}_{\B{P}^1})$, and so it belongs to $\Delta_{0}^{\rm p}$, yielding
	\[
	\begin{array}{ccccl}
		H_0\cdot\delta_{0}^{\rm p}&=&\K{C}^{0}\cdot\delta_{0}&=&12
		\\[2mm]
		H_0\cdot\delta_{0}^{\rm b}&=&\K{C}^{0}\cdot\delta_{0}&=&12
	\end{array}
	\]
	All other intersection numbers are $0$, except for $F_0\cdot\lambda=1$, $G_0\cdot\lambda=H_0\cdot\lambda=3$. In summary, we obtain a table:
	\[
	\begin{array}{c|cccccccc}
		&\lambda&\delta_{0}^{\rm t}&\delta_{0}^{\rm p}&\delta_{0}^{\rm b}&\delta_{1}^{\rm n}&\delta_{1}^{\rm t}&\delta_{1}^{\rm p}&\delta_{(j\geq 2)}\\[\dimexpr-\normalbaselineskip+0.5mm]
		\\\hline
		\\[\dimexpr-\normalbaselineskip+1mm]
		F_0&1&0&12&0&0&-1&0&0\\
		G_0&3&12&0&12&-3&0&0&0\\
		H_0&3&0&12&12&0&0&-3&0
	\end{array}
	\]
	Note that the formulas
	\[
	\begin{array}{ccccl}
		\pi_{*}(G_0)\cdot\delta_{0}&=&G_0\cdot\pi^{*}\delta_{0}&=&G_0\cdot(\delta_{0}^{\rm t}+2\,\delta_{0}^{\rm b})
		\\[2mm]
		\pi_{*}(H_0)\cdot\delta_{0}&=&H_0\cdot\pi^{*}\delta_{0}&=&H_0\cdot(\delta_{0}^{\rm p}+2\,\delta_{0}^{\rm b})
	\end{array}
	\]
	both hold.
\end{exam}

\begin{exam}[irreducible nodal curves]\label{testirrednod}
	In keeping with the notation used in example \ref{BoundDivRgO}, we fix a general curve $(B,p)\in\K{M}_{g-1,\,1}$ and consider the test curve obtained by gluing $p$ to a varying point $y\in B$, namely
	\[
	\K{Y}=\N{Bl}_{(p,p)}(B\times B)\slash(\Delta_{B}\sim B\times\{p\})\lto B
	\]
	This corresponds to a family
	\[
	\begin{array}{rclclcl}
		\K{Y}&\equiv&\{B_{py}\}_{y\in B}&\subset&\Delta_{0}&\subset&\OK{M}_{g}
	\end{array}
	\]
	where $B_{py}$ is an irreducible nodal curve for $y\neq p$ and $B_{pp}$ is a copy of $B$ with a \textsl{pigtail} attached to $p$, in the sense of \cite{HMModuliCurves} Section 3.C. Again, we can readily see that the intersection table of $\K{Y}$ with the generators of $\Pic(\OK{M}_g)_{\B{Q}}$ is:
	\[
	\begin{array}{c|cccccccc}
		&\lambda&\delta_{0}&\delta_{1}&\delta_{(j\geq 2)}\\[\dimexpr-\normalbaselineskip+0.5mm]
		\\\hline
		\\[\dimexpr-\normalbaselineskip+1mm]
		\K{Y}&0&2-2g&1&0
	\end{array}
	\]
	Pulling back $\K{Y}$ by the map $\Delta_{0}^{\rm t}\to\Delta_{0}$, we lift it to a test curve $Y_{0}$ such that:
	\[
	\begin{array}{rclclcl}
		Y_0&\equiv&\{(B_{py},\,\eta_{y}^{\rm t})\;\slash\;\eta_{y}^{\rm t}\in\Delta_{0}^{\rm t}(B_{py}) \}_{y\in B}&\subset&\Delta_{0}^{\rm t}&\subset&\OK{R}_g
	\end{array}
	\]
	Since $\deg(\Delta_{0}^{\rm t}\vert\Delta_{0})=1$, we have $\pi_{*}(Y_0)=\K{Y}$, hence
	\[
	\begin{array}{ccccccc}
		Y_0\cdot \delta_{0}^{\rm t}&=&Y_0\cdot\pi^{*}\delta_{0}&=&\K{Y}\cdot\delta_{0}&=&2-2g
	\end{array}
	\]
	In addition, the special fiber $\eta_{p}^{\rm t}$ lies in $\Delta_{1}^{\rm n}$, as it pulls back to the trivial bundle $(\K{O}_{B},\K{O}_{\B{P}^{1}})$ on the normalization $B\times\B{P}^{1}$ of $B_{pp}$, and thus is trivial over $B$. Then the last non-zero intersection number standing is
	\[
	\begin{array}{ccccccc}
		Y_0\cdot \delta_{1}^{\rm n}&=&Y_0\cdot\pi^{*}\delta_{1}&=&\K{Y}\cdot\delta_{1}&=&1
	\end{array}
	\]
	and we get a table:
	\[
	\begin{array}{c|cccccccc}
		&\lambda&\delta_{0}^{\rm t}&\delta_{0}^{\rm p}&\delta_{0}^{\rm b}&\delta_{1}^{\rm n}&\delta_{1}^{\rm t}&\delta_{1}^{\rm p}&\delta_{(j\geq 2)}\\[\dimexpr-\normalbaselineskip+0.5mm]
		\\\hline
		\\[\dimexpr-\normalbaselineskip+1mm]
		Y_0&0&2-2g&0&0&1&0&0&0
	\end{array}
	\]
	Note that we could have also pulled back by $\Delta_{0}^{\rm p}\to\Delta_{0}$ or $\Delta_{0}^{\rm b}\to\Delta_{0}$.
\end{exam}

%% file: 220Dnull.tex
\input{221IntroDnull}

\subsection[The Prym-null divisor and its irreducible components]{The divisor $\K{P}_{\mathbf{null}}$ and its irreducible components}
\label{SubSectPrymNull}

\input{222DefDnull}

\subsection{Parity change: smooth case}
\label{SubSectPrymSmooth}

\input{223Parity}

\subsection{Parity change: irreducible nodal case}
\label{SubSectPrymIrred}

\input{224Paritynodal}

%% file: 221IntroDnull.tex
In this section, we define the even and odd Prym-null divisors and study how a theta characteristic changes parity when tensored by a Prym root.

%% file: 222DefDnull.tex
Let $C$ be a smooth, integral curve of genus $g$.

\begin{defn}\label{vanishingthetanull}
	An even theta characteristic $\theta$ on $C$ with $h^{0}(C,\theta)\neq 0$ (that is, with $h^{0}(C,\theta)\geq 2$ and $h^{0}(C,\theta)\equiv 0\,$ mod $2$) is called a \textsl{vanishing theta-null}.
\end{defn}

The terminology here may seem confusing, as vanishing theta-nulls are even theta characteristics with \textsl{non-vanishing} global sections. This is justified by the classical theory of theta functions, whose \textsl{Thetanullwert} vanishes only when the associated even theta characteristic is a vanishing theta-null; see \cite{BeauvilleThetanull}.

The locus of curves with a vanishing theta-null, namely
\[
\Theta_{\rm null}=\{(C,\theta)\in\K{S}_g^+\;\slash\;h^0(C,\theta)\geq 2\}=\K{S}_g^+\cap\K{W}_{g-1,\,g}^2
\]
gives rise to the \textsl{theta-null divisor} $\Theta_{\rm null}$ on $\K{S}_{g}^{+}$, as well as its closure $\OC{\Theta}_{\rm null}$ in $\OK{S}_{g}^{+}$. This divisor plays an important role in the study of the geometry of $\OK{S}_g^+$, due to its effective nature and geometric characterization: for example, a computation of the class of $\OC{\Theta}_{\rm null}$ allows \cite{FarkasEvenSpin} to prove that $\OK{S}_{g}^{+}$ is of general type for $g\geq 9$, and of non-negative Kodaira dimension if $g=8$.

The theta-null divisor can be pushed forward by $\pi_{+}\colon\K{S}_{g}^{+}\to\K{M}_{g}$ to obtain
\[
\K{M}_g^{\rm null}=\{C\in\K{M}_g\;\slash\;\exists\,\theta\in S_g^+(C)\textrm{ with }h^0(C,\theta)\geq 2\}\subset\K{M}_{g}
\]
whose closure $\OK{M}_g^{\rm null}$ in $\OK{M}_{g}$ is described by \cite{TeixidorMgnull}. In turn, pulling back $\K{M}_g^{\rm null}$ by $\pi_{\K{R}}\colon\K{R}_{g}\to\K{M}_{g}$ results in a divisor
\[
\K{P}_{\rm null}=\{(C,\eta)\in\K{R}_g\;\slash\;\exists\,\theta\in S_g^+(C)\textrm{ with }h^0(C,\theta)\geq 2\}\subset\K{R}_{g}
\]
Note that the line bundle $\theta\otimes\eta$ is again a theta characteristic, different from $\theta$, which may therefore be even or odd. Moreover, \cite{TeixidorHalfCan} shows that the projection $\Theta_{\rm null}\to\K{M}_{g}^{\rm null}$ is generically finite of degree $1$, hence we can build a rational map
\[
\K{P}_{\rm null}\to\K{S}_{g}=\K{S}_{g}^{+}\sqcup\K{S}_{g}^{-},\qquad(C,\eta)\mapsto(C,\theta\otimes\eta)
\]
where $\theta\in\Theta_{\rm null}(C)$. Then, with the temporary notation
\[
\bar{\theta}=\theta\otimes\eta\in S_{g}(C),\qquad\theta=\bar{\theta}\otimes\eta\in\Theta_{\rm null}(C)
\]
we may rewrite the defining condition of $\K{P}_{\rm null}$ as
\[
\K{P}_{\rm null}=\{(C,\eta)\in\K{R}_g\;\slash\;\exists\,\bar{\theta}\in S_g(C)\textrm{ with }\bar{\theta}\otimes\eta\in\Theta_{\rm null}(C)\}\subset\K{R}_{g}
\]
and deduce that the parity of $\bar{\theta}=\theta\otimes\eta$ yields a decomposition
\[
\K{P}_{\rm null}=\K{P}_{\rm\! null}^{+}\sqcup\K{P}_{\rm\! null}^{-}
\]
Dropping the bar for the sake of simplicity, we get the following:

\begin{defn}\label{defDnullevenodd}
	We refer to the divisor $\K{P}_{\rm null}$ on $\K{R}_{g}$ as the \textsl{Prym-null divisor}. Accordingly, its irreducible components $\K{P}_{\rm\! null}^{+}$ and $\K{P}_{\rm\! null}^{-}$, namely
	\[
	\begin{array}{c}
		\K{P}_{\rm\! null}^{+}=\{(C,\eta)\in\K{R}_g\;\slash\;\exists\,\theta\in S_g^+(C)\textrm{ with }\theta\otimes\eta\in\Theta_{\rm null}(C)\}\subset\K{R}_{g}
		\\[2mm]
		\K{P}_{\rm\! null}^{-}=\{(C,\eta)\in\K{R}_g\;\slash\;\exists\,\theta\in S_g^-(C)\textrm{ with }\theta\otimes\eta\in\Theta_{\rm null}(C)\}\subset\K{R}_{g}
	\end{array}
	\] 
	with $\K{P}_{\rm null}=\K{P}_{\rm\! null}^{+}+\K{P}_{\rm\! null}^{-}$, are called the \textsl{even} and \textsl{odd Prym-null divisors}.
\end{defn}

Since the Prym-null divisors are natural, geometric divisors on $\K{R}_{g}$, our goal is to compute the class of their closures $\OK{P}_{\rm\! null}^{+}$, $\OK{P}_{\rm\! null}^{-}$ in $\OK{R}_{g}$. Such a computation would continue the work of \cite{TeixidorMgnull} and \cite{FarkasEvenSpin}, where the classes of $\OK{M}_{g}^{\rm null}$ and $\OC{\Theta}_{\rm null}$ are respectively expressed in terms of the generating classes of $\Pic(\OK{M}_{g})_{\B{Q}}$ and $\Pic(\OK{S}_{g}^{+})_{\B{Q}}$. In particular, write
\[
\OC{\mu}_g^{\rm null}=\K{O}_{\OK{M}_{g}}(\OK{M}_g^{\rm null})\in\Pic(\OK{M}_{g}),\qquad\OC{\vartheta}_{\rm null}=\K{O}_{\OK{S}_{g}^{+}}(\OC{\Theta}_{\rm null})\in\Pic(\OK{S}_{g}^{+}).
\]
for the aforementioned classes, and consider the notation $\lambda,\delta_{i},\delta_{i}^{\rm x}$ introduced in examples \ref{BoundDivSgI} and \ref{BoundDivSgO}. Then \cite{TeixidorMgnull} and \cite{FarkasEvenSpin} provide formulas
\[
\begin{array}{rcl}
	\OC{\mu}_g^{\rm null}&=&\displaystyle 2^{g-3}\,\Big((2^g+1)\,\lambda-2^{g-3}\,\delta_{0}-\sum_{i=1}^{\left\lfloor g/2 \right\rfloor}(2^i-1)(2^{g-i}-1)\,\delta_i\Big)
	\\[1mm]
	\OC{\vartheta}_{\rm null}&=&\displaystyle\frac{1}{4}\,\lambda-\frac{1}{16}\,\delta_{0}^{\rm n}-\frac{1}{2}\,\sum_{i=1}^{\left\lfloor g/2 \right\rfloor}\delta_{i}^{-}
\end{array}
\]
the latter of which implies the former, as the class $[\OK{M}_{g}^{\rm null}]$ can also be obtained by pushing forward the class $[\OC{\Theta}_{\rm null}]$ by the coarse moduli map $\OC{S}_g^+\to\OC{M}_g$.

Let us write the classes of $\OK{P}_{\rm\! null}$, $\OK{P}_{\rm\! null}^{+}$ and $\OK{P}_{\rm\! null}^{-}$ as
\[
\begin{array}{rcl}
	\varrho_{\rm null}=\K{O}_{\OK{R}_{g}}(\OK{P}_{\rm\! null})\in\Pic(\OK{R}_{g}),&\quad&\varrho_{\rm null}=\varrho^{+}_{\rm null}+\varrho^{-}_{\rm null}
	\\[1mm]
	\varrho^{+}_{\rm null}=\K{O}_{\OK{R}_{g}}(\OK{P}_{\rm\! null}^{+})\in\Pic(\OK{R}_{g})\phantom{,}
	\\[1mm]
	\varrho^{-}_{\rm null}=\K{O}_{\OK{R}_{g}}(\OK{P}_{\rm\! null}^{-})\in\Pic(\OK{R}_{g})\phantom{,}
\end{array}
\]
and recall the notation $\delta_{0}^{\rm t},\delta_{0}^{\rm p},\delta_{0}^{\rm b},\delta_{i}^{\rm n},\delta_{i}^{\rm t},\delta_{i}^{\rm p}$ from examples \ref{BoundDivRgI} and \ref{BoundDivRgO}. The sum $\varrho_{\rm null}$ can be directly computed as the pullback of $\OC{\mu}_{g}^{\rm null}$ by the map
\[
\begin{array}{rcrcl}
	\pi_{\K{R}}\colon\OK{R}_{g}\to\OK{M}_{g},
	&\quad&\pi_{\K{R}}^{*}(\lambda)&\hspace*{-4pt}=\hspace*{-4pt}&\lambda
	\\[1mm]
	&&\pi_{\K{R}}^{*}(\delta_{0})&\hspace*{-4pt}=\hspace*{-4pt}&\delta_{0}^{\rm t}+\delta_{0}^{\rm p}+2\,\delta_{0}^{\rm b}
	\\[1mm]
	&&\pi_{\K{R}}^{*}(\delta_{i})&\hspace*{-4pt}=\hspace*{-4pt}&\delta_{i}^{\rm n}+\delta_{i}^{\rm t}+\delta_{i}^{\rm p}
\end{array}
\]
with $1\leq i<g/2$, and moreover $\pi_{\K{R}}^{*}(\delta_{g/2})=\delta_{g/2}^{\rm n}+\delta_{g/2}^{\rm p}$ for even $g$.

\begin{prop}\label{classprymnull}
	The class of $\OK{P}_{\rm\! null}$ in $\Pic(\OK{R}_{g})_{\B{Q}}$ is given by
	\[
	\begin{array}{rcccl}
		\varrho_{\rm null}
		&=&
		2^{g-3}\,\Big((2^g+1)\,\lambda&-&2^{g-3}\,(\delta_{0}^{\rm t}+\delta_{0}^{\rm p}+2\,\delta_{0}^{\rm b})
		\vspace*{2mm}\\
		&&&-&\displaystyle\sum_{i=1}^{k}(2^i-1)(2^{g-i}-1)(\delta_{i}^{\rm n}+\delta_{i}^{\rm t}+\delta_{i}^{\rm p})
		\vspace*{2mm}\\
		&&&-&\psi(g)\cdot (2^{g/2}-1)^{2}(\delta_{g/2}^{\rm n}+\delta_{g/2}^{\rm p})\Big)
	\end{array}
	\]
	where the upper bound $k$ and the parity-checking function $\psi(g)$, defined as
	\[
	\begin{array}{rclcl}
		k&=&\left\lceil g/2 \right\rceil-1&=&
		\left\lbrace
		\begin{array}{ll}
			\left\lfloor g/2 \right\rfloor&\textrm{\rm{ if $g$ odd}}\\[1mm]
			\left\lfloor g/2 \right\rfloor-1&\textrm{\rm{ if $g$ even}}
		\end{array}\right.
		\\[6mm]
		\psi(g)&=&\dfrac{1+(-1)^{g}}{2}&=&
		\left\lbrace
		\begin{array}{ll}
			\mathrlap{0}\phantom{\left\lfloor g/2 \right\rfloor-1}&\textrm{\rm{ if $g$ odd}}	\\[1mm]
			1&\textrm{\rm{ if $g$ even}}
		\end{array}\right.
	\end{array}
	\]
	account for the slight variation in pullback that occurs when $g=2i$.
\end{prop}
\begin{dem}
	Follows from $\pi_{\K{R}}^{*}(\OC{\mu}_{g}^{\rm null})=\varrho_{\rm null}$ and the formulas above. \qed
\end{dem}

\begin{rem}\label{doublecheck}
	Once the classes $\varrho^{+}_{\rm null}$ and $\varrho^{-}_{\rm null}$ are computed, proposition \ref{classprymnull} offers a quick double-check of their accuracy, by virtue of $\varrho^{+}_{\rm null}+\varrho^{-}_{\rm null}=\varrho_{\rm null}$.
\end{rem}

\begin{rem}\label{defPnullcoeffs}
	With the notation of proposition \ref{classprymnull}, we may write
	\[
	\begin{array}{rcccl}
		\varrho^{+}_{\rm null}
		&=&
		\lambda^{+}\cdot\lambda&-&\left(\delta_{0}^{\rm t, +}\cdot\delta_{0}^{\rm t}+\delta_{0}^{\rm p, +}\cdot\delta_{0}^{\rm p}+\delta_{0}^{\rm b, +}\cdot\delta_{0}^{\rm b}\right)
		\vspace*{2mm}\\
		&&&-&\displaystyle\sum_{i=1}^{k}(\delta_{i}^{\rm n, +}\cdot\delta_{i}^{\rm n}+\delta_{i}^{\rm t, +}\cdot\delta_{i}^{\rm t}+\delta_{i}^{\rm p, +}\cdot\delta_{i}^{\rm p})
		\vspace*{2mm}\\
		&&&-&\psi(g)\cdot(\delta_{g/2}^{\rm n, +}\cdot\delta_{g/2}^{\rm n}+\delta_{g/2}^{\rm p, +}\cdot\delta_{g/2}^{\rm p})
		\vspace*{5mm}\\
		\varrho^{-}_{\rm null}
		&=&
		\lambda^{-}\cdot\lambda&-&\left(\delta_{0}^{\rm t, -}\cdot\delta_{0}^{\rm t}+\delta_{0}^{\rm p, -}\cdot\delta_{0}^{\rm p}+\delta_{0}^{\rm b, -}\cdot\delta_{0}^{\rm b}\right)
		\vspace*{2mm}\\
		&&&-&\displaystyle\sum_{i=1}^{k}(\delta_{i}^{\rm n, -}\cdot\delta_{i}^{\rm n}+\delta_{i}^{\rm t, -}\cdot\delta_{i}^{\rm t}+\delta_{i}^{\rm p, -}\cdot\delta_{i}^{\rm p})
		\vspace*{2mm}\\
		&&&-&\psi(g)\cdot(\delta_{g/2}^{\rm n, -}\cdot\delta_{g/2}^{\rm n}+\delta_{g/2}^{\rm p, -}\cdot\delta_{g/2}^{\rm p})
	\end{array}
	\]
	and subsequently aim our efforts at determining the rational coefficients
	\[
	\begin{array}{rclcl}
		\lambda^{+},\,\delta_{0}^{\rm t, +},\,\delta_{0}^{\rm p, +},\,\delta_{0}^{\rm b, +},\,\delta_{i}^{\rm n, +},\,\delta_{i}^{\rm t, +},\,\delta_{i}^{\rm p, +}&\in&\B{Q}&\quad&\textrm{(resp. $-$)}
	\end{array}
	\]
	for $1\leq i\leq \left\lfloor g/2 \right\rfloor$. To that end, the assortment of test curves introduced earlier will prove to be most useful.
\end{rem}

According to definition \ref{defDnullevenodd}, the essential distinction between $\OK{P}_{\rm\! null}^{+}$ and $\OK{P}_{\rm\! null}^{-}$ is the parity of the theta characteristic $\theta\in S_{g}(C)$, whereas its associated tensored version $\theta\otimes\eta\in\Theta_{\rm null}(C)$ is always even. This leads towards the question of how the parity of a theta characteristic changes when it is tensored by a Prym root.

%% file: 223Parity.tex
Given a Prym pair $(C,\eta)$ of genus $g$, we can consider a map
\[
S_{g}(C)\lto S_{g}(C),\qquad\theta\lmto\theta\otimes\eta
\]
and wonder how the parity of $\theta$ is affected by it.

\begin{defn}
	Let $(C,\eta)$ be a Prym pair of genus $g$. Consider the subsets
	\[
	\begin{array}{rcccl}
		S_{\eta}^{+,+}(C)&=&\{\theta\in S_{g}^{+}(C)\;\slash\;\theta\otimes\eta\in S_{g}^{+}(C)\}&\subset&S_{g}^{+}(C)
		\\[2mm]
		S_{\eta}^{+,-}(C)&=&\{\theta\in S_{g}^{+}(C)\;\slash\;\theta\otimes\eta\in S_{g}^{-}(C)\}&\subset&S_{g}^{+}(C)
		\\[2mm]
		S_{\eta}^{-,+}(C)&=&\{\theta\in S_{g}^{-}(C)\;\slash\;\theta\otimes\eta\in S_{g}^{+}(C)\}&\subset&S_{g}^{-}(C)
		\\[2mm]
		S_{\eta}^{-,-}(C)&=&\{\theta\in S_{g}^{-}(C)\;\slash\;\theta\otimes\eta\in S_{g}^{-}(C)\}&\subset&S_{g}^{-}(C)
	\end{array}
	\]
	into which $S_{g}(C)$ decomposes as a disjoint union.
\end{defn}

\begin{rem}
	Note that $S_{\eta}^{+,-}(C)\cong S_{\eta}^{-,+}(C)$, $\theta\mapsto\theta\otimes\eta$. This leaves us with three distinct sets that we want to study.
\end{rem}

For any smooth, integral, genus $g$ curve $C$, the group $J_{2}(C)$ acts on $S_{g}(C)$ by means of the map
\[
J_2(C)\times S_{g}(C)\lto S_{g}(C),\qquad(\eta,\theta)\lmto\theta\otimes\eta
\]
whose associated difference map can be written as
\[
\N{diff}\colon\; S_{g}(C)\times S_{g}(C)\lto J_2(C),\qquad(\theta_1,\theta_2)\lmto\theta_1\otimes\theta_2^{-1}
\]
If we remove the diagonal $\Delta=\N{diff}^{-1}(\K{O}_X)$, we get a map
\[
\N{diff}_{\neq}\colon\; S_{g}(C)\times S_{g}(C)-\Delta\lto R_g(C)
\]
whose fibers, of order $2^{2g}$, reflect how many ways there are of writing a Prym root $\eta$ as a difference of theta characteristics $\theta_1\otimes\theta_2^{-1}$, that is, with $\theta_1=\theta_2\otimes\eta$. Since we aim to keep track of the parity of $\theta=\theta_{2}$ and $\theta\otimes\eta=\theta_{1}$, we just need to consider the restrictions
\[
\begin{array}{rlcl}
	\N{diff}_{+}\colon&S_{g}^{+}(C)\times S_{g}^{+}(C)-\Delta&\lto&R_g(C)
	\\[2mm]
	\N{diff}_{-}\colon&S_{g}^{-}(C)\times S_{g}^{-}(C)-\Delta&\lto&R_g(C)
	\\[2mm]
	\N{diff}_{\pm}\colon&S_{g}^{+}(C)\times S_{g}^{-}(C)&\lto&R_g(C)
\end{array}
\]
Recalling that
\[
\begin{array}{lcl}
	\#R_{g}(C)=2^{2g}-1,&\quad&\#S_{g}^{+}(C)=2^{g-1}(2^{g}+1)
	\\[1mm]
	&\quad&\#S_{g}^{-}(C)=2^{g-1}(2^{g}-1)
\end{array}
\]
it is easy to count the fibers of these difference maps.

\begin{lem}
	With the previous notation, it holds that
	\[
	\begin{array}{rcl}
		\#\N{diff}_{+}^{-1}(\eta)&=&2^{g-1}(2^{g-1}+1)
		\\[2mm]
		\#\N{diff}_{-}^{-1}(\eta)&=&2^{g-1}(2^{g-1}-1)
		\\[2mm]
		\#\N{diff}_{\pm}^{-1}(\eta)&=&2^{2g-2}
	\end{array}
	\]
	for any Prym pair $(C,\eta)$ of genus $g$.
\end{lem}
\begin{dem}
	These numbers follow from the computation
	\[\begin{array}{rcccl}
		\#\N{diff}_{+}^{-1}(\eta)&=&\dfrac{\#S_{g}^{+}(C)\cdot(\#S_{g}^{+}(C)-1)}{\#R_{g}(C)}&=&2^{g-1}(2^{g-1}+1)
		\\[6mm]
		\#\N{diff}_{-}^{-1}(\eta)&=&\dfrac{\#S_{g}^{-}(C)\cdot(\#S_{g}^{-}(C)-1)}{\#R_{g}(C)}&=&2^{g-1}(2^{g-1}-1)
		\\[6mm]
		\#\N{diff}_{\pm}^{-1}(\eta)&=&\dfrac{\#S_{g}^{+}(C)\cdot\#S_{g}^{-}(C)}{\#R_{g}(C)}&=&2^{2g-2}
	\end{array}
	\]
	which depends only on the genus $g$ of the curve. \qed
\end{dem}

\begin{defn}
	We denote
	\[
	\begin{array}{lclcl}
		\N{N}_{g}^{+}&=&\#\N{diff}_{+}^{-1}(\eta)&=&2^{g-1}(2^{g-1}+1)
		\\[2mm]
		\N{N}_{g}^{-}&=&\#\N{diff}_{-}^{-1}(\eta)&=&2^{g-1}(2^{g-1}-1)
		\\[2mm]
		\N{N}_{g}^{\pm}&=&\#\N{diff}_{\pm}^{-1}(\eta)&=&2^{2g-2}
	\end{array}
	\]
	for any positive integer $g\in\B{Z}^{+}$.
\end{defn}

\begin{prop}\label{Prymrootcount}
	Let $(C,\eta)$ be a Prym pair of genus $g$. Under the map
	\[
	S_{g}(C)\to S_{g}(C),\qquad\theta\mapsto\theta\otimes\eta
	\]
	that is, when tensoring by $\eta$, there are:
	\begin{enumerate}[label=\rm(\roman*)]
		\item $\N{N}_{g}^{+}=2^{g-1}(2^{g-1}+1)$ even theta characteristics on $C$ that remain even.
		\item $\N{N}_{g}^{-}=2^{g-1}(2^{g-1}-1)$ odd theta characteristics on $C$ that remain odd.
		\item $\N{N}_{g}^{\pm}=2^{2g-2}$ even theta characteristics on $C$ that become odd.
		\item $\N{N}_{g}^{\pm}=2^{2g-2}$ odd theta characteristics on $C$ that become even.
	\end{enumerate}
	In particular, $\#S_{g}(C)=\N{N}_{g}^{+}+\N{N}_{g}^{-}+2\,\N{N}_{g}^{\pm}=2^{2g}$.
\end{prop}
\begin{dem}
	As suggested above, we have
	\[
	\begin{array}{cclcl}
		\#S_{\eta}^{+,+}(C)&=&\multicolumn{3}{l}{\#\{\theta\in S_{g}^{+}(C)\;\slash\;\theta\otimes\eta\in S_{g}^{+}(C)\}}
		\\[1mm]
		&=&\multicolumn{3}{l}{\#\{(\theta_1,\,\theta_2)\in S_{g}^{+}(C)\times S_{g}^{+}(C)\;\slash\;\theta_1=\theta_2\otimes\eta\in S_{g}^{+}(C)\}}
		\\[1mm]
		&=&\multicolumn{3}{l}{\#\{(\theta_1,\,\theta_2)\in S_{g}^{+}(C)\times S_{g}^{+}(C)-\Delta\;\slash\;\theta_1\otimes\theta_2^{-1}=\eta\}}
		\\[1mm]
		&=&\#\N{diff}_{+}^{-1}(\eta)&=&\N{N}_{g}^{+}
	\end{array}
	\]
	and similarly $\#S_{\eta}^{-,-}(C)=\N{N}_{g}^{-}$ and $\#S_{\eta}^{+,-}(C)=\#S_{\eta}^{-,+}(C)=\N{N}_{g}^{\pm}$. \qed
\end{dem}

%% file: 224Paritynodal.tex
Let $(B,p,q)\in\K{M}_{g-1,\,2}$ and take the irreducible nodal curve $X=B_{pq}\in\OK{M}_{g}$ obtained from $B$ by gluing the points $p$ and $q$, with normalization $\nu\colon B\to B_{pq}$. The dualizing bundle $\omega_{X}$ is the subbundle of $\nu_{*}(\omega_{B}(p+q))$ fulfilling the residue condition, that is, such that the following diagram commutes:
\[
\xymatrix{
&\omega_{X}\ar@{}[d]|>(.55){\circlearrowleft}
\ar[dl]_-{\N{res}_{p}}
\ar[dr]^-{\N{res}_{q}}\\
\kappa(p)
\ar@{}@<4pt>[rr]|-{\cong}
\ar[rr]_-{-1}&&
\kappa(q)
}
\]
In this particular case, we actually have $H^{0}(X,\omega_{X})=H^{0}(B,\omega_{B}(p+q))$, since
\[
\begin{array}{rclclcl}
	h^{0}(\omega_{B}(p+q))&=&2g-2-(g-1)+1&=&g&=&h^{0}(\omega_{X})
\end{array}
\]
by Riemann-Roch and duality. As mentioned in example \ref{BoundDivSgO}, a spin curve
\[
(X,\theta,\alpha)\in\OK{S}_{g},\qquad\alpha\colon\theta^{\otimes 2}\cong\omega_{X},\qquad(\nu^{*}\theta)^{\otimes 2}\cong\omega_{B}(p+q)
\]
is given by a root $\theta_{B}\in\sqrt{\omega_{B}(p+q)}$ and a suitable gluing $\varphi\colon\restr{\theta_{B}}{p}\cong\restr{\theta_{B}}{q}$, which by the above discussion is bound to a condition $\varphi^{\otimes 2}\equiv -1$ corresponding to the commutativity of the following diagram:
\[
\xymatrix{
	\restr{\omega_{B}(p+q)}{p}
	\ar[rr]^-{\varphi^{\otimes 2}}_-{\cong}
	\ar[d]_-{\N{res}_{p}}^-{\cong}&
	\ar@{}[d]|-{\circlearrowleft}&
	\restr{\omega_{B}(p+q)}{q}
	\ar[d]^-{\N{res}_{q}}_-{\cong}\\
	\kappa(p)
	\ar@{}@<4pt>[rr]|-{\cong}
	\ar[rr]_-{-1}&&
	\kappa(q)
}
\]
Specifically, consider the canonical isomorphism $\psi$ induced by the diagram
\[
\xymatrix@R=8pt{
&&&\restr{\theta_{B}}{p}
\ar[dd]^-{\psi}_{\cong}\ar[dr]\\
0\ar[r]&\theta_{B}(-p-q)\ar[r]
&\theta_{B}\ar@{}[r]|-{\circlearrowright}
\ar[ur]\ar[dr]&&0\\
&&&\restr{\theta_{B}}{q}\ar[ur]
}
\]
where $\restr{\theta_{B}}{p}$ and $\restr{\theta_{B}}{q}$ are expressed as cokernels of the same map. Let us give an explicit description of $\psi$. On the one hand, Riemann-Roch and duality yield
\[
\begin{array}{rclcl}
	h^{0}(\theta_{B})-h^{0}(\theta_{B}(-p-q))&=&g-1-(g-1)+1&=&1
\end{array}
\]
so we can write $\restr{\theta_{B}}{p}=\langle\sigma(p)\rangle$ and $\restr{\theta_{B}}{q}=\langle\sigma(q)\rangle$ for any section
\[
\sigma\;\in\; H^{0}(B,\theta_{B})-H^{0}(B,\theta_{B}(-p-q))
\]
and see that $\psi$ is, by definition, the isomorphism
\[
\psi\colon\restr{\theta_{B}}{p}\cong\restr{\theta_{B}}{q}\,,\qquad\sigma(p)\mapsto\sigma(q)
\]
On the other hand, we have $\sigma^{\otimes 2}\in H^{0}(B,\omega_{B}(p+q))=H^{0}(X,\omega_{X})$, hence
\[
\N{res}_{p}(\sigma^{\otimes 2})+\N{res}_{q}(\sigma^{\otimes 2})=0
\]
and it is clear that $\psi^{\otimes 2}\equiv -1$, in the sense of:
\[
\xymatrix@C=16pt{
	\restr{\omega_{B}(p+q)}{p}
	\ar[rr]^-{\psi^{\otimes 2}}_-{\cong}
	\ar[d]_-{\N{res}_{p}}^-{\cong}&
	\ar@{}[d]|-{\circlearrowleft}&
	\restr{\omega_{B}(p+q)}{q}
	\ar[d]^-{\N{res}_{q}}_-{\cong}
	&&
	\sigma^{\otimes 2}(p)
	\ar@{|->}[rr]
	\ar@{|->}[d]&
	\ar@{}[d]|-{\circlearrowleft}&
	\sigma^{\otimes 2}(q)
	\ar@{|->}[d]
	\\
	\kappa(p)
	\ar@{}@<4pt>[rr]|-{\cong}
	\ar[rr]_-{-1}&&
	\kappa(q)
	&&
	\N{res}_{p}(\sigma^{\otimes 2})
	\ar@{|->}[rr]&&
	\N{res}_{q}(\sigma^{\otimes 2})
}
\]
If we also consider the opposite isomorphism
\[
-\psi\colon\restr{\theta_{B}}{p}\cong\restr{\theta_{B}}{q}\,,\qquad\sigma(p)\mapsto-\sigma(q)\,,\qquad(-\psi)^{\otimes 2}\equiv -1
\]
then $\psi$ and $-\psi$ are the only possible ways of gluing $p$ and $q$ to make $\theta_{B}$ into a square root of $\omega_{X}$. The resulting bundles on $X$, which we denote by
\[
(\theta_{B},\psi),\;(\theta_{B},-\psi)\;\in\;\sqrt{\omega_{X}\vphantom{(p+q)}}
\]
are thus the two elements in the fiber $(\nu^{*})^{-1}(\theta_{B})$ of the double cover
\[
\nu^{*}\colon\sqrt{\omega_{X}\vphantom{(p+q)}}\lto\sqrt{\omega_{B}(p+q)}
\]
Furthermore, observe that the $1$-dimensional space of sections $\langle\sigma\rangle\subset H^{0}(B,\theta_{B})$ is preserved under the gluing $\psi$, but lost under the gluing $-\psi$. As a result, the dimension of the glued global sections is given by
\[
\begin{array}{rcl}
	h^{0}(X,(\theta_{B},\psi))&=&h^{0}(B,\theta_{B})
	\\[1mm]
	h^{0}(X,(\theta_{B},-\psi))&=&h^{0}(B,\theta_{B})-1
\end{array}
\]
so that $(\theta_{B},\psi)$ and $(\theta_{B},-\psi)$ always have different parity.

Finally, let $\eta^{\rm t}$ be the single Prym root of $X$ lying in the divisor $\Delta_{0}^{\rm t}\subset\OK{R}_{g}$, as defined in example \ref{BoundDivRgO}. In other words, we have
\[
\eta^{\rm t}\neq\K{O}_{X},\qquad\nu^{*}\eta^{t}=\K{O}_{B}
\]
It follows that tensoring by $\eta^{\rm t}$ permutes the elements of $(\nu^{*})^{-1}(\theta_{B})$, since
\[
(\theta_{B},\psi)\otimes\eta^{\rm t}\neq(\theta_{B},\psi),\qquad\nu^{*}((\theta_{B},\psi)\otimes\eta^{\rm t})=\theta_{B}
\]
and similarly for $(\theta_{B},-\psi)$. This corresponds to a change in parity:

\begin{prop}\label{irredparitychange}
	With the notation above, let $(X,\theta,\alpha)\in\OK{S}_{g}^{+}$ be a general point of $\Delta_{0}^{\rm n}$. Then tensoring by $(X,\eta^{\rm t},\beta)\in\Delta_{0}^{\rm t}\subset\OK{R}_{g}$ induces a new spin curve $(X,\theta\otimes\eta^{\rm t},\alpha\otimes\beta)$ in $\Delta_{0}^{\rm n}\subset\OK{S}_{g}^{-}$, of opposite parity (resp. $\OK{S}_{g}^{-}$, $\OK{S}_{g}^{+}$).
\end{prop}

%% file: 230ClassComp.tex
\input{231IntroClass}

\subsection{Over reducible nodal curves}
\label{SubSectPnullRed}

\input{232OverRedNod}

\subsection{Over curves with elliptic tails}
\label{SubSectPnullEllip}

\input{233OverEllipt}

\subsection{Over irreducible nodal curves}
\label{SubSectPnullIrred}

\input{234OverIrredNod}

\subsection{Class expansion and application to other families}
\label{SubSectPnullClass}

\input{235ClassDnull}

%% file: 231IntroClass.tex
In this section, we use test curve techniques to determine all the coefficients in the rational expansions of the Prym-null classes.

%% file: 232OverRedNod.tex
Recall the test curves from example \ref{testrednod}, that is
\[
\begin{array}{rccclcl}
	F_i&\equiv&\{(C\cup_{y\sim q}D,\,(\eta_C,\K{O}_D))\}_{y\in C}&\subset&\Delta_{i}^{\rm n}&\subset&\OK{R}_g
	\\[2mm]
	G_i&\equiv&\{(C\cup_{y\sim q}D,\,(\K{O}_C,\eta_D))\}_{y\in C}&\subset&\Delta_{i}^{\rm t}&\subset&\OK{R}_g
	\\[2mm]
	H_i&\equiv&\{(C\cup_{y\sim q}D,\,(\eta_C,\eta_D))\}_{y\in C}&\subset&\Delta_{i}^{\rm p}&\subset&\OK{R}_g
	\\[3mm]
	\multicolumn{2}{r}{\rm with}&
	\multicolumn{5}{l}{
	C\in\K{M}_i,\;(D,q)\in\K{M}_{g-i,\,1}\;\textrm{ general,}}
	\\[1mm]
	\multicolumn{2}{r}{\rm and}&
	\multicolumn{5}{l}{\eta_C\in R_i(C), 
		\;\eta_D\in R_{g-i}(D)\;\textrm{ arbitrary.}}
\end{array}
\]
If $g\neq 2i$, their intersection table is:
\[
\begin{array}{c|cccccccc}
	&\lambda&\delta_{0}^{\rm t}&\delta_{0}^{\rm p}&\delta_{0}^{\rm b}&\delta_{i}^{\rm n}&\delta_{i}^{\rm t}&\delta_{i}^{\rm p}&\delta_{(j\neq i)}\\[\dimexpr-\normalbaselineskip+0.5mm]
	\\\hline
	\\[\dimexpr-\normalbaselineskip+1mm]
	F_i&0&0&0&0&2-2i&0&0&0\\
	G_i&0&0&0&0&0&2-2i&0&0\\
	H_i&0&0&0&0&0&0&2-2i&0
\end{array}
\]
If $g=2i$, we have $F_i\cdot\delta_{i}^{\rm n}=G_i\cdot\delta_{i}^{\rm n}=2-2i$ instead.

We want to determine $F_i\cap\OK{P}_{\rm\! null}^{+}$.

\begin{rem}\label{LimitSmoothFam}
	If a stable Prym curve
	\[
	\begin{array}{rclcl}
		F_{i,\,y}&=&(C\cup_{y\sim q}D,\,(\eta_C,\K{O}_D))&\in&F_i
	\end{array}
	\]
	lies in $\OK{P}_{\rm\! null}^{+}$, then it can be expressed as the limit of a smooth family in $\K{P}_{\rm\! null}^{+}$, in the following sense. First, let us write
	\[
	\begin{array}{rclclcl}
		(X_{y},\eta_{y})&=&\multicolumn{5}{l}{(C\cup_{y\sim 0}E\cup_{q\sim \infty}D,\,(\eta_C,\K{O}_{E},\K{O}_D))}
		\\[2mm]
		\N{st}(X_{y},\eta_{y})&=&(C\cup_{y\sim q}D,\,(\eta_C,\K{O}_D))&=&F_{i,\,y}&\in&\OK{P}_{\rm\! null}^{+}
	\end{array}
	\]
	to account for the exceptional component $E\cong\B{P}^{1}$ which appears when working with stable spin structures on $C\cup_{y\sim q}D$. Then there exist families
	\[
	\begin{array}{rrclcl}
		\multicolumn{4}{c}{
			f\colon\K{X}\to \Spec(R)=\{\xi,\,p_0\}}
		&&\textrm{of quasistable curves,}
		\\[2mm]
		\;\;&(\N{st}(f),\,\eta,\,\beta)&\hspace*{-6pt}\in\hspace*{-6pt}&\OK{R}_{g}&&\textrm{of stable Prym curves, and}
		\\[2mm]
		&(f,\,\theta,\,\alpha)&\hspace*{-6pt}\in\hspace*{-6pt}&\OK{S}_{g}^{+}&&\textrm{of stable (even) spin curves,}
	\end{array}
	\]
	such that:
	\begin{enumerate}[label=\rm(\roman*)]
		\item $\K{X}$ is a smooth surface.
		\item $R$ is a discrete valuation ring with maximal ideal $\G{m}$, whose spectrum is composed of a special point $p_{0}\equiv\G{m}$ and a generic point $\xi\equiv(0)$.
		\item On the special fiber $\K{X}_{0}=f^{-1}(p_0)$, it holds that
		\[
		\begin{array}{rcl}
			(\K{X}_{0},\restr{\eta}{\K{X}_{0}})=(X_{y},\eta_{y}),&
			\qquad&
			\N{st}(\K{X}_{0},\restr{\eta}{\K{X}_{0}})=F_{i,\,y}\in\OK{P}_{\rm\! null}^{+}
		\end{array}
		\]
		\item On the generic fiber $\K{X}_{\xi}=f^{-1}(\xi)=\N{st}(f)^{-1}(\xi)$, it holds that
		\[
		\begin{array}{rcl}
			(\K{X}_{\xi},\eta_{\xi})\in\K{P}_{\rm\! null}^{+},&
			\qquad&
			(\K{X}_{\xi},\theta_{\xi}\otimes\eta_{\xi})\in\Theta_{\rm null}
		\end{array}
		\]
		or equivalently $(\theta_{\xi}\otimes\eta_{\xi})^{\otimes 2}\simeq\omega_{\K{X}_{\xi}}$ and $h^{0}(\K{X}_{\xi},\theta_{\xi}\otimes\eta_{\xi})\geq 2$, $\,\equiv 0\,$ mod $2$.
	\end{enumerate}
	If we recall example \ref{BoundDivSgI} and use the notation
	\[
	\begin{array}{rclclcl}
		\restr{\theta}{\K{X}_{0}}&=&\theta_{y}^{+}&=&(\theta_C,\K{O}_E(1),\theta_D)&\in&\OK{S}_{g}^{+}(X_{y})
	\end{array}
	\]
	then it is clear that $\theta_C$ and $\theta_{D}$ must have the same parity. In addition, since $C$ and $D$ are general, the dimension of the global sections of $\theta_C$ and $\theta_{D}$ is at most one, and thus we get
	\[
	\begin{array}{rclcl}
		h^{0}(C,\theta_C)&=&h^{0}(D,\theta_{D})&\in&\{0,1\}
	\end{array}
	\]
	Observe that if $F_{i,\,y}$ were to lie in $\OK{P}_{\rm\! null}^{-}$, then $\restr{\theta}{\K{X}_{0}}=\theta_{y}^{-}$ would be odd instead of even and these theta characteristics would have opposite parity.
\end{rem}

As described in \cite{EisenHarrisLimit} Section 2, the data given in remark \ref{LimitSmoothFam} produces a limit linear series of type $\G{g}_{g-1}^{1}$ on $C\cup_{y\sim q}D$, namely
\[
\begin{array}{rclcl}
	\ell&=&\Big(\ell_C=(L_C,V_C),\;\ell_D=(L_D,V_D)\Big)&\in&G_{g-1}^{1}(C)\times G_{g-1}^{1}(D)
\end{array}
\]
where the line bundles $L_{C}$ and $L_{D}$ are obtained by looking at the equality
\[
\begin{array}{rclcl}
	\restr{\theta}{\K{X}_{0}}\otimes\restr{\eta}{\K{X}_{0}}&=&\theta_{y}^{+}\otimes\eta_y&=&(\theta_C\otimes\eta_C,\K{O}_{E}(1),\theta_D)
\end{array}
\]
and twisting with $y$ and $q$ to adjust the degrees to $g-1$, so that
\[
\left\lbrace
\begin{array}{lcl}
	L_C&=&\theta_C\otimes\eta_C\,((g-i)y)\vspace*{2mm}\\
	L_D&=&\theta_D\,(iq)
\end{array}
\right.
\]
Since $\theta_{\xi}\otimes\eta_{\xi}\in\Theta_{\rm null}(\K{X}_{\xi})$ is even and parity is constant in families, we get
\[
\begin{array}{rclclcl}
	h^{0}(\theta_C\otimes\eta_C)+h^{0}(\theta_D)&=&h^{0}(\theta_{y}^{+}\otimes\eta_y)&\equiv&h^{0}(\theta_{\xi}\otimes\eta_{\xi})&\equiv&0 \mod 2
\end{array}
\]
In particular, $\theta_C\otimes\eta_C$ and $\theta_D$ must have the same parity, and the dimension of their global sections is again either $0$ or $1$ due to generality. This results in two distinct possibilities for the $\OK{P}_{\rm\! null}^{+}$ setting, and two more for the $\OK{P}_{\rm\! null}^{-}$ one:
\[
\begin{array}{cccccccc}
	&
	&h^{0}(\theta_C)&h^{0}(\theta_D)&h^{0}(\theta_C\otimes\eta_C)
	\\[\dimexpr-\normalbaselineskip+0.5mm]
	\\\cline{3-5}
	&\\[\dimexpr-\normalbaselineskip+1mm]
	\multirow{2}{*}{$\OK{P}_{\rm\! null}^{+}$}
	&&0&0&0&&\rightsquigarrow&\N{(F,+,0)}\\
	&&1&1&1&&\rightsquigarrow&\N{(F,+,1)}
	\\[\dimexpr-\normalbaselineskip+0mm]
	&\\\cline{3-5}
	&\\[\dimexpr-\normalbaselineskip+1mm]
	\multirow{2}{*}{$\OK{P}_{\rm\! null}^{-}$}
	&&1&0&0&&\rightsquigarrow&\N{(F,-,0)}\\
	&&0&1&1&&\rightsquigarrow&\N{(F,-,1)}
	\\[\dimexpr-\normalbaselineskip+0mm]
	&\\\cline{3-5}
\end{array}
\]

In order to study each of these cases, we first need to recall a basic property of the vanishing sequence of a linear series.

\begin{rem}
	Given a linear series $(L,V)\in G^{r}_{d}(C)$ on a smooth curve $C$ and a point $p\in C$, we can find an ordered basis
	\[
	V=\langle s_{0},\ldots,s_{r}\rangle\subset H^{0}(C,L)
	\]
	such that, if we write $a_{i}(p)=\N{ord}_{p}(s_{i})$ for all $i\in\{0,\ldots,r\}$, then
	\[
	a_{0}(p)<\ldots<a_{r}(p)
	\]
	is the vanishing sequence of $(L,V)$ at $p$. Taking any $b\in\B{Z}^{+}$ and observing that $V(-bp)$ is the subspace of sections $s\in V$ such that $\N{ord}_{p}(s)\geq b$, we can extract a basis of $V(-bp)$ out of $\langle s_{0},\ldots,s_{r}\rangle$, namely
	\[
	V(-bp)=\langle s_{j},\ldots,s_{r}\rangle\subset H^{0}(C,L(-bp))
	\]
	where the index $j\in\{0,\ldots,r+1\}$ is determined by the inequalities
	\[
	a_{j}(p)=\N{ord}_{p}(s_{j})\geq b,\qquad a_{j-1}(p)=\N{ord}_{p}(s_{j-1})<b
	\]
	whenever they make sense. Finally, the fact that there are $(r+1)-j$ elements in such a basis leads to the useful relation
	\[
	\begin{array}{rcl}
		\dim V(-bp)=r+1-j&\Leftrightarrow&a_{j-1}(p)<b\leq a_{j}(p)
	\end{array}
	\]
	which we will systematically use in the subsequent discussion. For example, we can apply it to $L_C=\theta_C\otimes\eta_C\,((g-i)y)$ and deduce that
	\[
	\begin{array}{rcl}
		h^{0}(\theta_C\otimes\eta_C)=h^{0}(L_{C}(-(g-i)y)=2-j&\Leftrightarrow& a_{j-1}^{\ell_C}(y)<g-i\leq a_{j}^{\ell_C}(y)
	\end{array}
	\]
	with $j\in\{1,2\}$ depending on the parity of $\theta_C\otimes\eta_C$.
\end{rem}

Let us start by analysing the two possibilities related to the even Prym-null divisor, labelled as in the table above.

\begin{poss}[$\N{F,+,0}$]
	In this case, we have
	\[
	\begin{array}{rcccl}
		h^{0}(\theta_C\otimes\eta_C)=0&\Rightarrow&a_{0}^{\ell_C}(y)<a_{1}^{\ell_C}(y)<g-i&\Rightarrow&a_{0}^{\ell_C}(y)\leq g-i-2
		\\[2mm]
		h^{0}(\theta_D)=0&\Rightarrow&a_{1}^{\ell_D}(q)<i&\Rightarrow&a_{1}^{\ell_D}(q)\leq i-1
	\end{array}
	\]
	Combining these upper bounds, we immediately reach a contradiction with one of the limit $\G{g}_{g-1}^{1}$ compatibility conditions
	\[
	g-1\leq a_{0}^{\ell_C}(y)+a_{1}^{\ell_D}(q)\leq g-3\qquad(!!)
	\]
	which therefore prevents this type of intersection from taking place.
\end{poss}

\begin{poss}[$\N{F,+,1}$]
	In this case, we have
	\[
	\begin{array}{rcl}
		h^{0}(\theta_C\otimes\eta_C)=1&\Rightarrow&a_{0}^{\ell_C}(y)<g-i\leq a_{1}^{\ell_C}(y)
		\\[2mm]
		h^{0}(\theta_D)=1&\Rightarrow&a_{0}^{\ell_D}(q)<i\leq a_{1}^{\ell_D}(q)
	\end{array}
	\]
	On the one hand, $(D,q)$ is general, so we may assume that $q\notin\N{supp}(\theta_D)$. Then the vanishing sequence of $\ell_D$ at $q$ can be bounded further:
	\[
	\begin{array}{rcl}
		h^{0}(\theta_D\,(-q))=h^{0}(\theta_D)-1=0&\Rightarrow&a_{1}^{\ell_D}(q)<i+1
		\\[2mm]
		h^{0}(\theta_D\,(q))=h^{0}(\theta_D)=1&\Rightarrow&a_{0}^{\ell_D}(q)<i-1\leq a_{1}^{\ell_D}(q)
	\end{array}
	\]
	We thus get $a_{1}^{\ell_D}(q)=i$ and, by the limit $\G{g}_{g-1}^{1}$ condition, $a_{0}^{\ell_C}(y)=g-i-1$. On the other hand, $C$ is general, so $\N{supp}(\theta_C\otimes\eta_C)$ consists of $i-1$ distinct points. As a result, we obtain a tight upper bound for $a_{1}^{\ell_C}(y)$, namely
	\[
	\begin{array}{rcl}
		\N{div}(\theta_C\otimes\eta_C)\ngeq 2y&\Rightarrow&\N{div}(s)\ngeq(g-i+2)y\quad\forall\,s\in H^{0}(L_{C})
		\\[2mm]
		&\Rightarrow&a_{1}^{\ell_C}(y)\leq g-i+1
	\end{array}
	\]
	which together with the condition $a_{0}^{\ell_D}(q)+a_{1}^{\ell_C}(y)\geq g-1$ yields $a_{0}^{\ell_D}(q)=i-2$ and $a_{1}^{\ell_C}(y)=g-i+1$. In turn, this means that $y\in\N{supp}(\theta_C\otimes\eta_C)$, and that $\ell$ is a refined limit $\G{g}_{g-1}^{1}$ of the form
	\[
	\left\lbrace
	\begin{array}{cclcl}
		\ell_C&=&\vert\theta_C\otimes\eta_C\,(y)\vert+(g-i-1)y&\in&G_{g-1}^{1}(C)
		\vspace*{2mm}\\
		\ell_D&=&\vert\theta_D\,(2q)\vert+(i-2)q&\in&G_{g-1}^{1}(D)
	\end{array}
	\right.
	\]
	with vanishing sequences $(g-i-1,\,g-i+1)$ and $(i-2,\,i)$.
\end{poss}

In conclusion, for each pair of theta characteristics $\theta_C$, $\theta_D$ fulfilling $\N{(F,+,1)}$, that is, such that $\theta_C\in S_{i}^{-}(C)$, $\theta_D\in S_{g-i}^{-}(D)$ and $\theta_C\otimes\eta_C\in S_{i}^{-}(C)$, then every $y\in\N{supp}(\theta_C\otimes\eta_C)$ yields a limit $\G{g}_{g-1}^{1}$ as above, and these limit linear series are the only ones contributing to the intersection $F_i\cap\OK{P}_{\rm\! null}^{+}$. Consequently, we need to count such pairs of theta characteristics.

Fortunately, we already have all the necessary tools to do this.

\begin{lem}\label{internumbFi}
	For all $i\in\{2,\ldots,g-1\}$, it holds that
	\[
	\begin{array}{rcl}
		F_i\cdot\OK{P}_{\rm\! null}^{+}&=&2^{g-2}(2^{i-1}-1)(2^{g-i}-1)(i-1)
	\end{array}
	\]
\end{lem}
\begin{dem}
	In light of the previous considerations, we may split the count into three parts. Specifically, we want to compute the number of:
	\begin{enumerate}[label=\rm(\roman*)]
		\item \textsl{Theta characteristics $\theta_C\in S_{i}^{-}(C)$ such that $\theta_C\otimes\eta_C\in S_{i}^{-}(C)$.}\\
		According to proposition \ref{Prymrootcount}, this is $\N{N}_{i}^{-}=2^{i-1}(2^{i-1}-1)$.
		\item \textsl{Theta characteristics $\theta_D\in S_{g-i}^{-}(D)$.}\\
		This is $\#S_{g-i}^{-}(D)=2^{g-i-1}(2^{g-i}-1)$.
		\item \textsl{Once $\theta_{C}$ is fixed, points $y$ in the support of $\theta_C\otimes\eta_C$.}\\
		Since $\theta_C\otimes\eta_C\in S_{i}(C)$, there are $\deg(\theta_C\otimes\eta_C)=i-1$ such points.
	\end{enumerate}
	Altogether, we obtain
	\[
	\begin{array}{rcl}
		F_i\cdot\OK{P}_{\rm\! null}^{+}&=&\#
		\left\lbrace
		\begin{array}{c}
			(\theta_C,\,\theta_D,\,y)\in S_{i}^{-}(C)\times S_{g-i}^{-}(D)\times C\;\slash\\[1mm]
			\,\theta_C\otimes\eta_C\in S_{i}^{-}(C),\,y\in\N{supp}(\theta_C\otimes\eta_C)
		\end{array}
		\right\rbrace\vspace*{2mm}\\
		&=&\N{N}_{i}^{-}\cdot\#S_{g-i}^{-}(D)\cdot(i-1)\\[2mm]
		&=&2^{i-1}(2^{i-1}-1)\cdot 2^{g-i-1}(2^{g-i}-1)\cdot(i-1)\\[2mm]
		&=&2^{g-2}(2^{i-1}-1)(2^{g-i}-1)(i-1)
	\end{array}
	\]
	as stated above. \qed
\end{dem}

In order to determine $F_i\cap\OK{P}_{\rm\! null}^{-}$ we follow the same argument, with the only difference being that $\theta_C$ and $\theta_D$ now have opposite parity (remark \ref{LimitSmoothFam}). Since this brings about minimal variations, we merely outline the situation and carry out the corresponding count. There are again two possibilities to tackle.

\begin{poss}[$\N{F,-,0}$]
	Similar contradiction to that of $\N{(F,+,0)}$.
\end{poss}

\begin{poss}[$\N{F,-,1}$]
	As with its even counterpart, we are able to build a limit linear series contributing to $F_i\cdot\OK{P}_{\rm\! null}^{-}$ whenever $y\in\N{supp}(\theta_C\otimes\eta_C)$. In this case, however, we have $\theta_C\in S_{i}^{+}(C)$.
\end{poss}

We thus get
\[
\begin{array}{rcl}
	F_i\cdot\OK{P}_{\rm\! null}^{-}&=&\#
	\left\lbrace
	\begin{array}{c}
		(\theta_C,\,\theta_D,\,y)\in S_{i}^{+}(C)\times S_{g-i}^{-}(D)\times C\;\slash\vspace*{1mm}\\
		\,\theta_C\otimes\eta_C\in S_{i}^{-}(C),\,y\in\N{supp}(\theta_C\otimes\eta_C)
	\end{array}
	\right\rbrace\vspace*{2mm}\\
	&=&\N{N}_{i}^{\pm}\cdot\#S_{g-i}^{-}(D)\cdot(i-1)\vspace*{2mm}\\
	&=&2^{2i-2}\cdot 2^{g-i-1}(2^{g-i}-1)\cdot(i-1)\vspace*{2mm}\\
	&=&2^{g+i-3}(2^{g-i}-1)(i-1)
\end{array}
\]

The procedure we have employed to study the intersection of $F_{i}$ and each of the Prym-null divisors works with $G_{i}$ and $H_{i}$ as well. Nonetheless, we still need to carefully track the small changes that happen along the way.

Let us briefly do this. If a stable Prym curve
\[
\begin{array}{rclcl}
	G_{i,\,y}&=&(C\cup_{y\sim q}D,\,(\K{O}_C,\eta_D))&\in&G_i
\end{array}
\]
lies in $\OK{P}_{\rm\! null}^{+}$ (resp. $\OK{P}_{\rm\! null}^{-}$), we can produce a limit $\G{g}_{g-1}^{1}$ on $C\cup_{y\sim q}D$ such that
\[
\left\lbrace
\begin{array}{lcl}
	L_C&=&\theta_C\,((g-i)y)\\[2mm]
	L_D&=&\theta_D\otimes\eta_D\,(iq)
\end{array}
\right.
\]
with $h^{0}(\theta_C)+h^{0}(\theta_D\otimes\eta_D)\equiv 0 \mod 2$, where $\theta_C$ and $\theta_D$ have the same parity (resp. opposite parity). Then $\theta_C$ and $\theta_D\otimes\eta_D$ have the same parity and we get the following possibilities:
\[
\begin{array}{cccccccccl}
	&
	&h^{0}(\theta_C)&h^{0}(\theta_D)&h^{0}(\theta_D\otimes\eta_D)
	\\[\dimexpr-\normalbaselineskip+0.5mm]
	\\\cline{3-5}
	&\\[\dimexpr-\normalbaselineskip+1mm]
	\multirow{2}{*}{$\OK{P}_{\rm\! null}^{+}$}
	&&0&0&0&&\rightsquigarrow&\N{(G,+,0)}&
	\hspace*{-2.5mm}:\hspace*{-2mm}&\textrm{contradiction}\\
	&&1&1&1&&\rightsquigarrow&\N{(G,+,1)}&
	\hspace*{-2.5mm}:\hspace*{-2mm}& y\in\N{supp}(\theta_C)
	\\[\dimexpr-\normalbaselineskip+0mm]
	&\\\cline{3-5}
	&\\[\dimexpr-\normalbaselineskip+1mm]
	\multirow{2}{*}{$\OK{P}_{\rm\! null}^{-}$}
	&&0&1&0&&\rightsquigarrow&\N{(G,-,0)}&
	\hspace*{-2.5mm}:\hspace*{-2mm}&\textrm{contradiction}\\
	&&1&0&1&&\rightsquigarrow&\N{(G,-,1)}&
	\hspace*{-2.5mm}:\hspace*{-2mm}& y\in\N{supp}(\theta_C)
	\\[\dimexpr-\normalbaselineskip+0mm]
	&\\\cline{3-5}
\end{array}
\]

With the only contribution of $\N{(G,+,1)}$ and $\N{(G,-,1)}$ to their respective intersections, we obtain
\[
\begin{array}{rcl}
	G_i\cdot\OK{P}_{\rm\! null}^{+}&=&\#
	\left\lbrace
	\begin{array}{c}
		(\theta_C,\,\theta_D,\,y)\in S_{i}^{-}(C)\times S_{g-i}^{-}(D)\times C\;\slash\\[1mm]
		\,\theta_D\otimes\eta_D\in S_{g-i}^{-}(D),\,y\in\N{supp}(\theta_C)
	\end{array}
	\right\rbrace\vspace*{2mm}\\
	&=&\#S_{i}^{-}(C)\cdot\N{N}_{g-i}^{-}\cdot(i-1)\\[2mm]
	&=&2^{i-1}(2^{i}-1)\cdot 2^{g-i-1}(2^{g-i-1}-1)\cdot(i-1)\\[2mm]
	&=&2^{g-2}(2^{i}-1)(2^{g-i-1}-1)(i-1)
\end{array}
\]
and similarly
\[\begin{array}{rcl}
	G_i\cdot\OK{P}_{\rm\! null}^{-}&=&\#
	\left\lbrace
	\begin{array}{c}
		(\theta_C,\,\theta_D,\,y)\in S_{i}^{-}(C)\times S_{g-i}^{+}(D)\times C\;\slash\\[1mm]
		\,\theta_D\otimes\eta_D\in S_{g-i}^{-}(D),\,y\in\N{supp}(\theta_C)
	\end{array}
	\right\rbrace\vspace*{2mm}\\
	&=&\#S_{i}^{-}(C)\cdot\N{N}_{g-i}^{\pm}\cdot(i-1)\\[2mm]
	&=&2^{i-1}(2^{i}-1)\cdot 2^{2g-2i-2}\cdot(i-1)\\[2mm]
	&=&2^{2g-i-3}(2^{i}-1)(i-1)
\end{array}
\]

\begin{rem}\label{contradict}
	The contradiction in $\N{(G,+,0)}$ and $\N{(G,-,0)}$ is again
	\[
	g-1\leq a_{0}^{\ell_C}(y)+a_{1}^{\ell_D}(q)\leq g-3\qquad(!!)
	\]
	In general, this condition will fail every time we try to use theta characteristics without global sections to build a limit $\G{g}_{g-1}^{1}$ on our reducible nodal curve, so in the future we will refrain from detailing it any further.
\end{rem}

Finally, if a stable Prym curve
\[
\begin{array}{rclcl}
	H_{i,\,y}&=&(C\cup_{y\sim q}D,\,(\eta_C,\eta_D))&\in&H_i
\end{array}
\]
lies in $\OK{P}_{\rm\! null}^{+}$ (resp. $\OK{P}_{\rm\! null}^{-}$), we can produce a limit $\G{g}_{g-1}^{1}$ on $C\cup_{y\sim q}D$ such that
\[
\left\lbrace
\begin{array}{lcl}
	L_C&=&\theta_C\otimes\eta_C\,((g-i)y)\\[2mm]
	L_D&=&\theta_D\otimes\eta_D\,(iq)
\end{array}
\right.
\]
with $h^{0}(\theta_C\otimes\eta_C)+h^{0}(\theta_D\otimes\eta_D)\equiv 0 \mod 2$, where $\theta_C$ and $\theta_D$ have the same parity (resp. opposite parity). Then $\theta_C\otimes\eta_C$ and $\theta_D\otimes\eta_D$ have the same parity and we get the following possibilities:
\[
\begin{array}{ccccccccl}
	&
	&h^{0}(\theta_C)&h^{0}(\theta_D)&h^{0}(\theta_C\otimes\eta_C)&h^{0}(\theta_D\otimes\eta_D)
	\\[\dimexpr-\normalbaselineskip+0.5mm]
	\\\cline{3-6}
	&\\[\dimexpr-\normalbaselineskip+1mm]
	\multirow{4}{*}{$\OK{P}_{\rm\! null}^{+}$}
	&&0&0&0&0&&\rightsquigarrow&\textrm{contradiction}\\
	&&1&1&0&0&&\rightsquigarrow&\textrm{contradiction}\\
	&&0&0&1&1&&\rightsquigarrow&y\in\N{supp}(\theta_C\otimes\eta_C)\\
	&&1&1&1&1&&\rightsquigarrow&y\in\N{supp}(\theta_C\otimes\eta_C)
	\\[\dimexpr-\normalbaselineskip+0mm]
	&\\\cline{3-6}
	&\\[\dimexpr-\normalbaselineskip+1mm]
	\multirow{4}{*}{$\OK{P}_{\rm\! null}^{-}$}
	&&0&1&0&0&&\rightsquigarrow&\textrm{contradiction}\\
	&&1&0&0&0&&\rightsquigarrow&\textrm{contradiction}\\
	&&0&1&1&1&&\rightsquigarrow&y\in\N{supp}(\theta_C\otimes\eta_C)\\
	&&1&0&1&1&&\rightsquigarrow&y\in\N{supp}(\theta_C\otimes\eta_C)
	\\[\dimexpr-\normalbaselineskip+0mm]
	&\\\cline{3-6}
\end{array}
\]

The count has now grown in complexity, but not by much. We have
\[
\begin{array}{rcl}
	H_i\cdot\OK{P}_{\rm\! null}^{+}&=&\#
	\left\lbrace
	\begin{array}{c}
		(\theta_C,\,\theta_D,\,y)\in(S_{i}^{+}(C)\times S_{g-i}^{+}(D)\times C)\;\cup\hspace*{7mm}\\[1mm]
		\phantom{(\theta_C,\,\theta_D,\,y)\in}\hspace*{4mm}
		\cup\,(S_{i}^{-}(C)\times S_{g-i}^{-}(D)\times C)\;\slash
		\\[1mm]
		\,\theta_C\otimes\eta_C\in S_{i}^{-}(C),\,\theta_D\otimes\eta_D\in S_{g-i}^{-}(D),\\[1mm]
		\,y\in\N{supp}(\theta_C\otimes\eta_C)
	\end{array}
	\right\rbrace\vspace*{2mm}\\
	&=&(\N{N}_{i}^{\pm}\N{N}_{g-i}^{\pm}+\N{N}_{i}^{-}\N{N}_{g-i}^{-})\cdot(i-1)\\[2mm]
	&=&(2^{2i-2}\cdot 2^{2g-2i-2}+2^{i-1}(2^{i-1}-1)\cdot 2^{g-i-1}(2^{g-i-1}-1))\cdot(i-1)\\[2mm]
	&=&2^{g-2}(2^{g-1}-2^{i-1}-2^{g-i-1}+1)(i-1)
\end{array}
\]
and similarly
\[
\begin{array}{rcl}
	H_i\cdot\OK{P}_{\rm\! null}^{-}&=&\#
	\left\lbrace
	\begin{array}{c}
		(\theta_C,\,\theta_D,\,y)\in(S_{i}^{+}(C)\times S_{g-i}^{-}(D)\times C)\;\cup\hspace*{7mm}\\[1mm]
		\phantom{(\theta_C,\,\theta_D,\,y)\in}\hspace*{4mm}
		\cup\,(S_{i}^{-}(C)\times S_{g-i}^{+}(D)\times C)\;\slash
		\\[1mm]
		\,\theta_C\otimes\eta_C\in S_{i}^{-}(C),\,\theta_D\otimes\eta_D\in S_{g-i}^{-}(D),\\[1mm]
		\,y\in\N{supp}(\theta_C\otimes\eta_C)
	\end{array}
	\right\rbrace\vspace*{2mm}\\
	&=&(\N{N}_{i}^{\pm}\N{N}_{g-i}^{-}+\N{N}_{i}^{-}\N{N}_{g-i}^{\pm})\cdot(i-1)\\[2mm]
	&=&(2^{2i-2}\cdot 2^{g-i-1}(2^{g-i-1}-1)+2^{i-1}(2^{i-1}-1)\cdot 2^{2g-2i-2})\cdot(i-1)\\[2mm]
	&=&2^{g-2}(2^{g-1}-2^{i-1}-2^{g-i-1})(i-1)
\end{array}
\]

To summarize, we compile the simplified expressions for all three collections of intersections into the following extension of lemma \ref{internumbFi}.

\begin{lem}\label{internumbFGHi}
	For all $i\in\{2,\ldots,g-1\}$, we have intersection numbers
	\[
	\begin{array}{c|cccccccc}
		&	&\OK{P}_{\rm\! null}^{+}
		&\;	&\OK{P}_{\rm\! null}^{-}
		\\[\dimexpr-\normalbaselineskip+0.5mm]
		\\\hline
		\\[\dimexpr-\normalbaselineskip+1mm]
		F_i	&&(2^{i-1}-1)(2^{g-i}-1)\,r
		&&2^{i-1}(2^{g-i}-1)\,r\\
		G_i	&&(2^{i}-1)(2^{g-i-1}-1)\,r
		&&(2^{i}-1)\,2^{g-i-1}\,r\\
		H_i	&&(2^{g-1}-2^{i-1}-2^{g-i-1}+1)\,r
		&&(2^{g-1}-2^{i-1}-2^{g-i-1})\,r
	\end{array}
	\]
	where $r=2^{g-2}(i-1)=-2^{g-3}(2-2i)$.
\end{lem}

\begin{rem}
	A quick computation shows that these numbers pass the check suggested by remark \ref{doublecheck}. Indeed, we can easily see that
	\[
	\begin{array}{cclcr}
		F_i\cdot\OK{P}_{\rm\! null}^{+}+F_i\cdot\OK{P}_{\rm\! null}^{-}&=&2^{g-2}(2^{i}-1)(2^{g-i}-1)(i-1)&=&F_i\cdot\OK{P}_{\rm\! null}\\[2mm]
		G_i\cdot\OK{P}_{\rm\! null}^{+}+G_i\cdot\OK{P}_{\rm\! null}^{-}&=&2^{g-2}(2^{i}-1)(2^{g-i}-1)(i-1)&=&G_i\cdot\OK{P}_{\rm\! null}\\[2mm]
		H_i\cdot\OK{P}_{\rm\! null}^{+}+H_i\cdot\OK{P}_{\rm\! null}^{-}&=&2^{g-2}(2^{i}-1)(2^{g-i}-1)(i-1)&=&H_i\cdot\OK{P}_{\rm\! null}
	\end{array}
	\]
	where we have used example \ref{testrednod} and proposition \ref{classprymnull}.
\end{rem}

If we combine lemma \ref{internumbFGHi} with the intersection table of example \ref{testrednod}, we can derive the first batch of coefficients from the formulas in remark \ref{defPnullcoeffs}.

\begin{prop}\label{dicoeffs}
	Fix integers $g\geq 5$ and $i\in\{1,\ldots,\left\lfloor g/2 \right\rfloor\}$. The generating classes $\delta_{i}^{\rm n},\,\delta_{i}^{\rm t},\,\delta_{i}^{\rm p}\in\Pic(\OK{R}_{g})_{\B{Q}}$ have coefficients
	\[
	\begin{array}{rc|crcl}
		&&\;\!&\delta_{i}^{\rm n, +}&=&2^{g-3}(2^{i-1}-1)(2^{g-i}-1)
		\\[2mm]
		\varrho^{+}_{\rm null}
		&&&\delta_{i}^{\rm t, +}&=&2^{g-3}(2^{i}-1)(2^{g-i-1}-1)
		\\[2mm]
		&&&\delta_{i}^{\rm p, +}&=&2^{g-3}(2^{g-1}-2^{i-1}-2^{g-i-1}+1)
		\\[-1mm]
		\multicolumn{1}{c}{}
		\\
		&&&\delta_{i}^{\rm n, -}&=&2^{g-3}\,2^{i-1}(2^{g-i}-1)
		\\[2mm]
		\varrho^{-}_{\rm null}
		&&&\delta_{i}^{\rm t, -}&=&2^{g-3}(2^{i}-1)\,2^{g-i-1}
		\\[2mm]
		&&&\delta_{i}^{\rm p, -}&=&2^{g-3}(2^{g-1}-2^{i-1}-2^{g-i-1})
	\end{array}
	\]
	in the rational expansions of the Prym-null classes in genus $g$.
\end{prop}
\begin{dem}
	Every family of test curves generates a linear relation between the coefficients of each expansion. Due to the simplicity of the intersection tables of $F_{i}$, $G_{i}$ and $H_{i}$ with the generators of $\Pic(\OK{R}_g)_{\B{Q}}$, their corresponding linear relations directly determine one coefficient (sometimes the same one, if two relations are linearly dependent). For the sake of brevity, we shall describe this computation simply in the $F_{i}$ case, as all others are analogous. To begin with, we have
	\[
	\begin{array}{rclcl}
		\deg\varrho^{+}_{\rm null}(F_{i})&=&F_{i}\cdot\OK{P}_{\rm\! null}^{+}&=&-2^{g-3}(2^{i-1}-1)(2^{g-i}-1)(2-2i)
	\end{array}
	\]
	by lemma \ref{internumbFGHi}. Furthermore, remark \ref{defPnullcoeffs} and example \ref{testrednod} show that
	\[
	\begin{array}{rclcl}
		\deg\varrho^{+}_{\rm null}(F_{i})&=&-\delta_{i}^{\rm n, +}\deg\delta_{i}^{\rm n}(F_{i})&=&-\delta_{i}^{\rm n, +}(2-2i)
	\end{array}
	\]
	with the convention $\delta_{i}^{\rm n, +}=\delta_{g-i}^{\rm t, +}$ for all $i$. Since only one coefficient survives, we can immediately extract it from the resulting equation:
	\[
	\begin{array}{rclcl}
		F_{i}&\;\fto\;&\delta_{i}^{\rm n, +}&=&2^{g-3}(2^{i-1}-1)(2^{g-i}-1)
	\end{array}
	\]
	For $2\leq i\leq\left\lfloor g/2 \right\rfloor$, each coefficient can be similarly computed by means of:
	\[
	\begin{array}{rclcl}
		F_{i}&\N{or}&G_{g-i}&\rightsquigarrow&
		\left\lbrace
		\begin{array}{lcl}
			\delta_{i}^{\rm n, +}&=&2^{g-3}(2^{i-1}-1)(2^{g-i}-1)
			\vspace*{2mm}\\
			\delta_{i}^{\rm n, -}&=&2^{g-3}\,2^{i-1}(2^{g-i}-1)
		\end{array}\right.\vspace*{3mm}\\
		G_{i}&\N{or}&F_{g-i}&\rightsquigarrow&
		\left\lbrace
		\begin{array}{lcl}
			\delta_{i}^{\rm t, +}&=&2^{g-3}(2^{i}-1)(2^{g-i-1}-1)
			\vspace*{2mm}\\
			\delta_{i}^{\rm t, -}&=&2^{g-3}(2^{i}-1)\,2^{g-i-1}
		\end{array}\right.\vspace*{3mm}\\
		H_{i}&\N{or}&H_{g-i}&\rightsquigarrow&
		\left\lbrace
		\begin{array}{lcl}
			\delta_{i}^{\rm p, +}&=&2^{g-3}(2^{g-1}-2^{i-1}-2^{g-i-1}+1)
			\vspace*{2mm}\\
			\delta_{i}^{\rm p, -}&=&2^{g-3}(2^{g-1}-2^{i-1}-2^{g-i-1})
		\end{array}\right.
	\end{array}
	\]
	For $i=1$, we do not have families $F_1$, $G_1$ and $H_1$. Nevertheless, we can use:
	\[
	\begin{array}{lcl}
		G_{g-1}&\rightsquigarrow&
		\left\lbrace
		\begin{array}{lcl}
			\delta_{1}^{\rm n, +}&=&0
			\vspace*{2mm}\\
			\delta_{1}^{\rm n, -}&=&2^{g-3}(2^{g-1}-1)
		\end{array}\right.\vspace*{3mm}\\
		F_{g-1}&\rightsquigarrow&
		\left\lbrace
		\begin{array}{lcl}
			\delta_{1}^{\rm t, +}&=&2^{g-3}(2^{g-2}-1)
			\vspace*{2mm}\\
			\delta_{1}^{\rm t, -}&=&2^{g-3}\,2^{g-2}
		\end{array}\right.\vspace*{3mm}\\
		H_{g-1}&\rightsquigarrow&
		\left\lbrace
		\begin{array}{lcl}
			\delta_{1}^{\rm p, +}&=&2^{g-3}\,2^{g-2}
			\vspace*{2mm}\\
			\delta_{1}^{\rm p, -}&=&2^{g-3}(2^{g-2}-1)
		\end{array}\right.
	\end{array}
	\]
	that is to say, the above formulas hold for $i=1$ as well. \qed
\end{dem}

\begin{rem}
	Observe that $\delta_{1}^{\rm n, +}=0$ is a consequence of the fact that genus $1$ curves do not have nontrivial odd theta characteristics, hence $\N{N}_{1}^{-}=0$.
\end{rem}

%% file: 233OverEllipt.tex
Recall the test curves from example \ref{testelliptail}, that is
\[
\begin{array}{rclclcl}
	F_0&\equiv&\{(C\cup_{p\sim\sigma(\lambda)}E_{\lambda},\,(\eta_C,\K{O}_{E_{\lambda}}))\}_{\lambda\in\B{P}^1}&\subset&\Delta_{1}^{\rm t}&\subset&\OK{R}_g
	\\[2mm]
	G_0&\equiv&\{(C\cup_{p\sim\sigma(\lambda)}E_{\lambda},\,(\K{O}_C,\eta_{E_{\lambda}}))\;\slash\;\eta_{E_{\lambda}}\in\gamma_1^{-1}(E_{\lambda})\}_{\lambda\in\B{P}^1}&\subset&\Delta_{1}^{\rm n}&\subset&\OK{R}_g
	\\[2mm]
	H_0&\equiv&\{(C\cup_{p\sim\sigma(\lambda)}E_{\lambda},\,(\eta_C,\eta_{E_{\lambda}}))\;\slash\;\eta_{E_{\lambda}}\in\gamma_1^{-1}(E_{\lambda})\}_{\lambda\in\B{P}^1}&\subset&\Delta_{1}^{\rm p}&\subset&\OK{R}_g
	\\[3mm]
	\multicolumn{2}{r}{\rm with}&
	\multicolumn{5}{l}{
		(C,p)\in\K{M}_{g-1,\,1}\;\textrm{ general,}}
	\\[1mm]
	&&
	\multicolumn{5}{l}{\{E_{\lambda}\}_{\lambda\in\B{P}^1}\subset\OK{M}_{1}\;\textrm{ general pencil of plane cubics with basepoint $\sigma$,}}
	\\[1mm]
	&&
	\multicolumn{5}{l}{\eta_C\in R_{g-1}(C)\;\textrm{ arbitrary,}}
	\\[1mm]
	\multicolumn{2}{r}{\rm and}&
	\multicolumn{5}{l}{\gamma_1\colon\OK{R}_{1,1}\to\OK{M}_{1,1}\;\textrm{ forgetful degree $3$ branched covering.}}
\end{array}
\]
Our goal is to expand their intersection table to include the Prym-null divisors, as we did in lemma \ref{internumbFGHi} for the previous families of test curves. In this case, it turns out that the new intersection numbers are much less imposing:
\[
\begin{array}{c|ccccccccc;{1pt/2pt}ccc}
	&\lambda&\delta_{0}^{\rm t}&\delta_{0}^{\rm p}&\delta_{0}^{\rm b}&\delta_{1}^{\rm n}&\delta_{1}^{\rm t}&\delta_{1}^{\rm p}&\delta_{(j\geq 2)}
	&\hspace*{-6pt}&\hspace*{-6pt}&\OK{P}_{\rm\! null}^{+}&\OK{P}_{\rm\! null}^{-}
	\\[\dimexpr-\normalbaselineskip+0.5mm]
	&&&&&&&&&\hspace*{-6pt}\\\hline
	&&&&&&&&&\hspace*{-6pt}
	\\[\dimexpr-\normalbaselineskip+1mm]
	F_0&1&0&12&0&0&-1&0&0
	&\hspace*{-6pt}&\hspace*{-6pt}&0&0\\
	G_0&3&12&0&12&-3&0&0&0
	&\hspace*{-6pt}&\hspace*{-6pt}&0&0\\
	H_0&3&0&12&12&0&0&-3&0
	&\hspace*{-6pt}&\hspace*{-6pt}&0&0
\end{array}
\]
Let us check the veracity of this claim.

First, we determine $F_0\cap\OK{P}_{\rm\! null}^{+}$ (resp. $\OK{P}_{\rm\! null}^{-}$). If a stable Prym curve
\[
\begin{array}{rclcl}
	F_{0,\,\lambda}&=&(C\cup_{p\sim e}E_{\lambda},\,(\eta_C,\K{O}_{E_{\lambda}}))&\in&F_0
\end{array}
\]
lies in $\OK{P}_{\rm\! null}^{+}$ (resp. $\OK{P}_{\rm\! null}^{-}$), we can produce a limit $\mathfrak{g}_{g-1}^{1}$ on $C\cup_{p\sim e}E_{\lambda}$ such that
\[
\left\lbrace
\begin{array}{lcl}
	L_C&=&\theta_C\otimes\eta_C\,(p)
	\\[2mm]
	L_{E_{\lambda}}&=&\theta_{E_{\lambda}}\,((g-1)e)
\end{array}
\right.
\]
with $h^{0}(\theta_C\otimes\eta_C)+h^{0}(\theta_{E_{\lambda}})\equiv 0 \mod 2$, where $\theta_C$ and $\theta_{E_{\lambda}}$ have the same parity (resp. opposite parity). Then $\theta_C\otimes\eta_C$ and $\theta_{E_{\lambda}}$ have the same parity and we get the following possibilities:
\[
\begin{array}{cccccccl}
	&
	&h^{0}(\theta_C)&h^{0}(\theta_{E_{\lambda}})&h^{0}(\theta_C\otimes\eta_C)
	\\[\dimexpr-\normalbaselineskip+0.5mm]
	\\\cline{2-5}
	&\\[\dimexpr-\normalbaselineskip+1mm]
	\multirow{2}{*}{$\OK{P}_{\rm\! null}^{+}$}
	&&0&0&0&&\rightsquigarrow&\textrm{contradiction}\\
	&&1&1&1&&\rightsquigarrow&\N{(F_0,+,1)}
	\\[\dimexpr-\normalbaselineskip+0mm]
	&\\\cline{2-5}
	&\\[\dimexpr-\normalbaselineskip+1mm]
	\multirow{2}{*}{$\OK{P}_{\rm\! null}^{-}$}
	&&1&0&0&&\rightsquigarrow&\textrm{contradiction}\\
	&&0&1&1&&\rightsquigarrow&\N{(F_0,-,1)}
	\\[\dimexpr-\normalbaselineskip+0mm]
	&\\\cline{2-5}
\end{array}
\]

Note that all four theta characteristics of a genus $1$ curve are of degree zero, hence the nontrivial ones have no global sections other than zero. In particular, the dimension of the global sections of $\theta_{E_{\lambda}}$ is given by
\[
\begin{array}{rcl}
	h^{0}(\theta_{E_{\lambda}})&=&
	\left\lbrace
	\begin{array}{cl}
		1&\textrm{if }\theta_{E_{\lambda}}=\K{O}_{E_{\lambda}}
		\\[1mm]
		0&\textrm{otherwise}
	\end{array}
	\right.
\end{array}
\]
for all $\lambda\in\B{P}^{1}$, meaning that the table above is comprehensive. Since half of the scenarios are already covered by remark \ref{contradict}, we just need to look at $\N{(F_0,+,1)}$ and $\N{(F_0,-,1)}$ to conclude. This can be done in one fell swoop.

\begin{posss}[$\N{F_0,+,1)}$, $\N{(F_0,-,1}$]
	In both of these cases, we have
	\[
	\begin{array}{rcl}
		h^{0}(\theta_C\otimes\eta_C)=1&\Rightarrow&a_{0}^{\ell_C}(p)<1\leq a_{1}^{\ell_C}(p)
		\\[2mm]
		h^{0}(\theta_{E_{\lambda}})=1&\Rightarrow&a_{0}^{\ell_{\smash{E_{\lambda}}}}(e)<g-1\leq a_{1}^{\ell_{\smash{E_{\lambda}}}}(e)
	\end{array}\]
	Now $(C,p)$ is general, so we may assume that $p\notin\N{supp}(\theta_C\otimes\eta_C)$. Therefore
	\[
	\begin{array}{rcl}
		h^{0}(\theta_C\otimes\eta_C\,(p))=h^{0}(\theta_C\otimes\eta_C)=1&\Rightarrow&a_{0}^{\ell_C}(p)<0\leq a_{1}^{\ell_C}(p)
	\end{array}\]
	Furthermore, $E_{\lambda}$ is a genus $1$ curve, so $\deg(\theta_{E_{\lambda}})=0$. As a result, we get
	\[
	\begin{array}{rcl}
		h^{0}(\theta_{E_{\lambda}}\,(-e))=0&\Rightarrow&a_{1}^{\ell_{\smash{E_{\lambda}}}}(e)<g
	\end{array}
	\]
	that is, $a_{1}^{\ell_{\smash{E_{\lambda}}}}(e)=g-1$, which contradicts the limit $\mathfrak{g}_{g-1}^{1}$ condition
	\[
	g-1\leq a_{0}^{\ell_C}(p)+a_{1}^{\ell_{\smash{E_{\lambda}}}}(e)< g-1\qquad(!!)
	\]
	(Alternatively, from $p\notin\N{supp}(\theta_C\otimes\eta_C)$ we may also deduce
	\[
	\begin{array}{rcl}
		h^{0}(\theta_C\otimes\eta_C\,(-p))=h^{0}(\theta_C\otimes\eta_C)-1=0&\Rightarrow&a_{1}^{\ell_C}(p)<2
	\end{array}
	\]
	which yields $a_{1}^{\ell_C}(p)=1$ and, via the limit $\mathfrak{g}_{g-1}^{1}$ condition, $a_{0}^{\ell_{\smash{E_{\lambda}}}}(e)=g-2$. This should prevent $\theta_{E_{\lambda}}$ from being trivial, due to the implication
	\[
	\begin{array}{rcl}
		h^{0}(\theta_{E_{\lambda}}\,(e))=h^{0}(\K{O}_{E_{\lambda}}\,(e))=1&\Rightarrow&a_{0}^{\ell_{\smash{E_{\lambda}}}}(e)<g-2\leq a_{1}^{\ell_{\smash{E_{\lambda}}}}(e)
	\end{array}
	\]
	However, the triviality of $\theta_{E_{\lambda}}$ is ensured by $h^{0}(\theta_{E_{\lambda}})=1$ and $\deg(\theta_{E_{\lambda}})=0$.)
\end{posss}

As every possibility leads to a contradiction, the intersections $F_0\cap\OK{P}_{\rm\! null}^{+}$ and $F_0\cap\OK{P}_{\rm\! null}^{-}$ are both empty, and thus $F_0\cdot\OK{P}_{\rm\! null}^{+}=F_0\cdot\OK{P}_{\rm\! null}^{-}=0$.

Next we study $G_0\cap\OK{P}_{\rm\! null}^{+}$ (resp. $\OK{P}_{\rm\! null}^{-}$). If a stable Prym curve
\[
\begin{array}{rclcl}
	G_{0,\,\lambda}&=&(C\cup_{p\sim e}E_{\lambda},\,(\K{O}_C,\eta_{E_{\lambda}}))&\in&G_0
\end{array}
\]
lies in $\OK{P}_{\rm\! null}^{+}$ (resp. $\OK{P}_{\rm\! null}^{-}$), we can produce a limit $\mathfrak{g}_{g-1}^{1}$ on $C\cup_{p\sim e}E_{\lambda}$ such that
\[
\left\lbrace
\begin{array}{lcl}
	L_C&=&\theta_C\,(p)
	\\[2mm]
	L_{E_{\lambda}}&=&\theta_{E_{\lambda}}\otimes\eta_{E_{\lambda}}\,((g-1)e)
\end{array}
\right.
\]
with $h^{0}(\theta_C)+h^{0}(\theta_{E_{\lambda}}\otimes\eta_{E_{\lambda}})\equiv 0 \mod 2$, where $\theta_C$ and $\theta_{E_{\lambda}}$ have the same parity (resp. opposite parity). Then $\theta_C$ and $\theta_{E_{\lambda}}\otimes\eta_{E_{\lambda}}$ have the same parity and we get the following possibilities:
\[
\begin{array}{cccccccl}
	&
	&h^{0}(\theta_C)&h^{0}(\theta_{E_{\lambda}})&h^{0}(\theta_{E_{\lambda}}\otimes\eta_{E_{\lambda}})
	\\[\dimexpr-\normalbaselineskip+0.5mm]
	\\\cline{2-5}
	&\\[\dimexpr-\normalbaselineskip+1mm]
	\multirow{2}{*}{$\OK{P}_{\rm\! null}^{+}$}
	&&0&0&0&&\rightsquigarrow&\textrm{contradiction}\\
	&&1&1&1&&\rightsquigarrow&\N{(G_0,+,1)}
	\\[\dimexpr-\normalbaselineskip+0mm]
	&\\\cline{2-5}
	&\\[\dimexpr-\normalbaselineskip+1mm]
	\multirow{2}{*}{$\OK{P}_{\rm\! null}^{-}$}
	&&0&1&0&&\rightsquigarrow&\textrm{contradiction}\\
	&&1&0&1&&\rightsquigarrow&\N{(G_0,-,1)}
	\\[\dimexpr-\normalbaselineskip+0mm]
	&\\\cline{2-5}
\end{array}
\]

Once more, remark \ref{contradict} addresses half of the table, and the remaining half is tackled similarly to, if not more easily than, its $F_{0}$ counterpart.

\begin{poss}[$\N{G_0,+,1}$]
	We may repeat $\N{(F_0,+,1)}$'s argument, but a simpler procedure is available in this case. Since $E_{\lambda}$ is a genus $1$ curve, it holds that
	\[
	\begin{array}{rcl}
		\left.
		\begin{array}{rcccl}
			h^{0}(\theta_{E_{\lambda}})&\hspace*{-5pt}=\hspace*{-5pt}&h^{0}(\theta_{E_{\lambda}}\otimes\eta_{E_{\lambda}})&\hspace*{-5pt}=\hspace*{-5pt}&1
			\\[2mm]
			\deg(\theta_{E_{\lambda}})&\hspace*{-5pt}=\hspace*{-5pt}&\deg(\theta_{E_{\lambda}}\otimes\eta_{E_{\lambda}})&\hspace*{-5pt}=\hspace*{-5pt}&0
		\end{array}
		\right\rbrace
		&\Rightarrow&
		\left\lbrace
		\begin{array}{rcrrcll}
			\theta_{E_{\lambda}}&\hspace*{-5pt}=\hspace*{-5pt}&\multicolumn{2}{r}{\theta_{E_{\lambda}}\otimes\eta_{E_{\lambda}}}&\hspace*{-5pt}=\hspace*{-5pt}&\K{O}_{E_{\lambda}}&\Rightarrow
			\\[2mm]
			&&\hspace*{11pt}\Rightarrow&\eta_{E_{\lambda}}&\hspace*{-5pt}=\hspace*{-5pt}&\K{O}_{E_{\lambda}}&\;(!!)
		\end{array}
		\right.
	\end{array}
	\]
	in direct contradiction with the nontriviality of a Prym root.
\end{poss}

\begin{poss}[$\N{G_0,-,1}$]
	Same contradiction as in $\N{(F_0,-,1)}$.
\end{poss}

Hence both intersections are empty again and $G_0\cdot\OK{P}_{\rm\! null}^{+}=G_0\cdot\OK{P}_{\rm\! null}^{-}=0$.

Finally, let us consider $H_0\cap\OK{P}_{\rm\! null}^{+}$ (resp. $\OK{P}_{\rm\! null}^{-}$). If a stable Prym curve
\[
\begin{array}{rclcl}
	H_{0,\,\lambda}&=&(C\cup_{p\sim e}E_{\lambda},\,(\eta_C,\eta_{E_{\lambda}}))&\in&H_0
\end{array}
\]
lies in $\OK{P}_{\rm\! null}^{+}$ (resp. $\OK{P}_{\rm\! null}^{-}$), we can produce a limit $\mathfrak{g}_{g-1}^{1}$ on $C\cup_{p\sim e}E_{\lambda}$ such that
\[
\left\lbrace
\begin{array}{lcl}
	L_C&=&\theta_C\otimes\eta_C\,(p)
	\\[2mm]
	L_{E_{\lambda}}&=&\theta_{E_{\lambda}}\otimes\eta_{E_{\lambda}}\,((g-1)e)
\end{array}
\right.
\]
with $h^{0}(\theta_C\otimes\eta_C)+h^{0}(\theta_{E_{\lambda}}\otimes\eta_{E_{\lambda}})\equiv 0 \mod 2$, where $\theta_C$ and $\theta_{E_{\lambda}}$ have the same parity (resp. opposite parity). Then $\theta_C\otimes\eta_C$ and $\theta_{E_{\lambda}}\otimes\eta_{E_{\lambda}}$ have the same parity and we get the following possibilities:
\[
\begin{array}{ccccccccl}
	&
	&h^{0}(\theta_C)&h^{0}(\theta_{E_{\lambda}})&h^{0}(\theta_C\otimes\eta_C)&h^{0}(\theta_{E_{\lambda}}\otimes\eta_{E_{\lambda}})
	\\[\dimexpr-\normalbaselineskip+0.5mm]
	\\\cline{3-6}
	&\\[\dimexpr-\normalbaselineskip+1mm]
	\multirow{4}{*}{$\OK{P}_{\rm\! null}^{+}$}
	&&0&0&0&0&&\rightsquigarrow&\textrm{contradiction}\\
	&&1&1&0&0&&\rightsquigarrow&\textrm{contradiction}\\
	&&0&0&1&1&&\rightsquigarrow&\N{(H_0,+,1,0)}\\
	&&1&1&1&1&&\rightsquigarrow&\N{(H_0,+,1,1)}
	\\[\dimexpr-\normalbaselineskip+0mm]
	&\\\cline{3-6}
	&\\[\dimexpr-\normalbaselineskip+1mm]
	\multirow{4}{*}{$\OK{P}_{\rm\! null}^{-}$}
	&&0&1&0&0&&\rightsquigarrow&\textrm{contradiction}\\
	&&1&0&0&0&&\rightsquigarrow&\textrm{contradiction}\\
	&&0&1&1&1&&\rightsquigarrow&\N{(H_0,-,1,1)}\\
	&&1&0&1&1&&\rightsquigarrow&\N{(H_0,-,1,0)}
	\\[\dimexpr-\normalbaselineskip+0mm]
	&\\\cline{3-6}
\end{array}
\]

Even though there are more cases, they are all covered by already discussed arguments. We point to those outside of the scope of remark \ref{contradict}.

\begin{posss}[$\N{H_0,+,1,0)}$, $\N{(H_0,-,1,0}$]
	Same contradiction as in $\N{(F_0,+,1)}$.
\end{posss}

\begin{posss}[$\N{H_0,+,1,1)}$, $\N{(H_0,-,1,1}$]
	Same contradiction as in $\N{(G_0,+,1)}$.
\end{posss}

As a result, it is clear that $H_0\cdot\OK{P}_{\rm\! null}^{+}=H_0\cdot\OK{P}_{\rm\! null}^{-}=0$ as well.

\begin{rem}
	The above computations, in combination with example \ref{testelliptail} and proposition \ref{classprymnull}, show that
	\[
	\begin{array}{rclcl}
		F_0\cdot\OK{P}_{\rm\! null}^{+}+F_0\cdot\OK{P}_{\rm\! null}^{-}&=&0&=&F_0\cdot\OK{P}_{\rm\! null}
	\end{array}
	\]
	(resp. $G_{0}$, $H_{0}$), as expected by remark \ref{doublecheck}.
\end{rem}

\begin{lem}\label{LIlinrel}
	In the setting of proposition \ref{dicoeffs}, the families $F_{0}$, $G_{0}$ and $H_{0}$ provide three linearly independent linear relations
	\[
	\begin{array}{lcl}
		F_{0}&\rightsquigarrow&
		\left\lbrace
		\begin{array}{lcl}
			\lambda^{+}-12\,\delta_{0}^{\rm p, +}&=&-2^{g-3}(2^{g-2}-1)
			\vspace*{2mm}\\
			\lambda^{-}-12\,\delta_{0}^{\rm p, -}&=&-2^{2g-5}
		\end{array}\right.\vspace*{3mm}\\
		G_{0}&\rightsquigarrow&
		\left\lbrace
		\begin{array}{lcl}
			\lambda^{+}-4\,\delta_{0}^{\rm t, +}-4\,\delta_{0}^{\rm b, +}&=&0
			\vspace*{2mm}\\
			\lambda^{-}-4\,\delta_{0}^{\rm t, -}-4\,\delta_{0}^{\rm b, -}&=&-2^{g-3}(2^{g-1}-1)
		\end{array}\right.\vspace*{3mm}\\
		H_{0}&\rightsquigarrow&
		\left\lbrace
		\begin{array}{lcl}
			\lambda^{+}-4\,\delta_{0}^{\rm p, +}-4\,\delta_{0}^{\rm b, +}&=&-2^{2g-5}
			\vspace*{2mm}\\
			\lambda^{+}-4\,\delta_{0}^{\rm p, -}-4\,\delta_{0}^{\rm b, -}&=&-2^{g-3}(2^{g-2}-1)
		\end{array}\right.
	\end{array}
	\]
	between the coefficients of $\lambda,\,\delta_{0}^{\rm t},\,\delta_{0}^{\rm p},\,\delta_{0}^{\rm b}\in\Pic(\OK{R}_{g})_{\B{Q}}$ in each expansion.
\end{lem}
\begin{dem}
	Follows from combining proposition \ref{dicoeffs} with the intersection of $F_{0}$ and the Prym-null divisors (resp. $G_{0}$, $H_{0}$). \qed
\end{dem}

%% file: 234OverIrredNod.tex
Recall the test curve from example \ref{testirrednod}, that is
\[\begin{array}{rclclcl}
	Y_0&\equiv&\{(B_{py},\,\eta_{y}^{\rm t})\;\slash\;\eta_{y}^{\rm t}\in\Delta_{0}^{\rm t}(B_{py}) \}_{y\in B}&\subset&\Delta_{0}^{\rm t}&\subset&\OK{R}_g
	\\[3mm]
	\multicolumn{2}{r}{\rm with}&
	\multicolumn{5}{l}{
		(B,p)\in\K{M}_{g-1,\,1}\;\textrm{ general,}}
	\\[1mm]
	&&
	\multicolumn{5}{l}{B_{py}=B\slash\{y\sim p\}\;\textrm{ irreducible nodal curve for $y\neq p$,}}
	\\[1mm]
	\multicolumn{2}{r}{\rm and}&
	\multicolumn{5}{l}{B_{pp}\;\textrm{ copy of $B$ with a pigtail attached to $p$.}}
\end{array}
\]
As we will see below, its extended intersection table is:
\[
\begin{array}{c|ccccccccc;{1pt/2pt}cccc}
	&\lambda&\delta_{0}^{\rm t}&\delta_{0}^{\rm p}&\delta_{0}^{\rm b}&\delta_{1}^{\rm n}&\delta_{1}^{\rm t}&\delta_{1}^{\rm p}&\delta_{(j\geq 2)}
	&\hspace*{-6pt}&\hspace*{-6pt}&\OK{P}_{\rm\! null}^{+}
	&\hspace*{-6pt}&\OK{P}_{\rm\! null}^{-}
	\\[\dimexpr-\normalbaselineskip+0.5mm]
	&&&&&&&&&\hspace*{-6pt}\\\hline
	&&&&&&&&&\hspace*{-6pt}
	\\[\dimexpr-\normalbaselineskip+1mm]
	Y_0&0&2-2g&0&0&1&0&0&0
	&\hspace*{-6pt}&\hspace*{-6pt}&0
	&\hspace*{-6pt}&2^{g-3}(2^{g-2}(g-3)+1)
\end{array}
\]
The reason that the Prym-null intersection gravitates entirely towards the odd side is the parity change explored in proposition \ref{irredparitychange}. With the help of both this result and the notation used to prove it, we can easily determine $Y_0\cap\OK{P}_{\rm\! null}^{+}$ and $Y_0\cap\OK{P}_{\rm\! null}^{-}$.

On the one hand, if a stable Prym curve
\[
\begin{array}{rclcl}
	Y_{0,\,y}&=&(B_{py},\,\eta_{y}^{\rm t})&\in&Y_0
\end{array}
\]
lies in the even Prym-null divisor $\OK{P}_{\rm\! null}^{+}$, then by definition there exists a stable vanishing theta-null $\theta_{y}\in\OC{\Theta}_{\rm null}\subset\OK{S}_{g}^{+}$ over $B_{py}$ such that $\theta_{y}\otimes\eta_{y}^{\rm t}\in\OK{S}_{g}^{+}$ is even. Since proposition \ref{irredparitychange} shows that $\theta_{y}$ and $\theta_{y}\otimes\eta_{y}^{\rm t}$ always have opposite parity, no such vanishing theta-null can exist, and therefore $Y_0\cdot\OK{P}_{\rm\! null}^{+}=0$.

On the other hand, if a stable Prym curve
\[
\begin{array}{rclcl}
	Y_{0,\,y}&=&(B_{py},\,\eta_{y}^{\rm t})&\in&Y_0
\end{array}
\]
lies in the odd Prym-null divisor $\OK{P}_{\rm\! null}^{-}$ instead, then by definition there exists a stable vanishing theta-null $\theta_{y}\in\OC{\Theta}_{\rm null}\subset\OK{S}_{g}^{+}$ over $B_{py}$ such that $\theta_{y}\otimes\eta_{y}^{\rm t}\in\OK{S}_{g}^{-}$ is odd. Proposition \ref{irredparitychange} now makes redundant the second part of the condition, so we just need to count how many $B_{py}$ admit a vanishing theta-null. This can be done by taking advantage of the formula for the theta-null class
\[
\begin{array}{rcl}
	\OC{\vartheta}_{\rm null}&=&\displaystyle\frac{1}{4}\,\lambda-\frac{1}{16}\,\delta_{0}^{\rm n}-\frac{1}{2}\,\sum_{i=1}^{\left\lfloor g/2 \right\rfloor}\delta_{i}^{-}
\end{array}
\]
introduced in section \ref{SectPrymPar} and originally given by \cite{FarkasEvenSpin}. If we consider the test curve $Y_0^{\rm n}$ obtained as the pullback of $\{B_{py}\}_{y\in B}$ by the divisor $\Delta_{0}^{\rm n}\subset\OK{S}_{g}^{+}$, i.e.
\[
\begin{array}{rclclcl}
	Y_0^{\rm n}&\equiv&\{(B_{py},\,\theta_{y})\;\slash\;\theta_{y}\in\Delta_{0}^{\rm n}(B_{py}) \}_{y\in B}&\subset&\Delta_{0}^{\rm n}&\subset&\OK{S}_g^{+}
\end{array}
\]
then the previous discussion identifies the intersection $Y_0\cap\OK{P}_{\rm\! null}^{-}$ with the intersection $Y_0^{\rm n}\cap\OC{\Theta}_{\rm null}$. The latter can easily be derived from the theta-null formula and the intersection table of $Y_0^{\rm n}$ with the generators of $\Pic(\OK{S}_g^{+})_{\B{Q}}$, which is:
\[
\begin{array}{c|cccccc}
	&\lambda&\delta_{0}^{\rm n}&\delta_{0}^{\rm b}&\delta_{1}^{+}&\delta_{1}^{-}&\delta_{(j\geq 2)}\\[\dimexpr-\normalbaselineskip+0.5mm]
	\\\hline
	\\[\dimexpr-\normalbaselineskip+1mm]
	Y_0^{\rm n}&0&2^{2g-1}(1-g)&0&2^{g-2}(2^{g-1}+1)&2^{g-2}(2^{g-1}-1)&0
\end{array}
\]
Indeed, example \ref{testirrednod} and $\deg(\Delta_{0}^{\rm n}\vert\Delta_{0})=2^{2g-2}$ yield all coefficients except for the $\delta_{1}^{+}$, $\delta_{1}^{-}$ ones, while over the special point $(B_{pp},\,\theta_{p})$ we can see that
\[
\begin{array}{rclcl}
	\theta_{p}&=&(\theta_{B},\,(\K{O}_{\B{P}^{1}},\varphi))&\in&\Delta_{0}^{\rm n}(B_{pp})
\end{array}
\]
for some $\theta_{B}\in S_{g-1}(B)$ and $\varphi\in\{\psi,-\psi\}$ such that $h^{0}(\theta_{B})\equiv h^{0}(\K{O}_{\B{P}^{1}},\varphi)\mod 2$. Since by construction of $\psi$ and $-\psi$ we have
\[
\begin{array}{rcl}
	h^{0}(\K{O}_{\B{P}^{1}},\varphi)&=&
	\left\lbrace
	\begin{array}{ccl}
		1&\Leftrightarrow&\varphi=\psi\\[1mm]
		0&\Leftrightarrow&\varphi=-\psi
	\end{array}
	\right.
\end{array}
\]
it follows that $\varphi$ is determined by the parity of $\theta_{B}$ and thus
\[
\begin{array}{rclcl}
	Y_0^{\rm n}\cdot\delta_{1}^{+}&=&\#S_{g-1}^{+}(B)&=&2^{g-2}(2^{g-1}+1)
	\\[2mm]
	Y_0^{\rm n}\cdot\delta_{1}^{-}&=&\#S_{g-1}^{-}(B)&=&2^{g-2}(2^{g-1}-1)
\end{array}
\]
In conclusion, we get
\[
\begin{array}{rclcl}
	Y_0\cdot\OK{P}_{\rm\! null}^{-}&=&Y_0^{\rm n}\cdot\OC{\Theta}_{\rm null}&=&-2^{-4}\,Y_0^{\rm n}\cdot\delta_{0}^{\rm n}-2^{-1}\,Y_0^{\rm n}\cdot\delta_{1}^{-}
	\\[2mm]
	&&&=&-2^{2g-5}(1-g)-2^{g-3}(2^{g-1}-1)
	\\[2mm]
	&&&=&2^{g-3}(2^{g-2}(g-3)+1)
\end{array}
\]
as indicated earlier.

\begin{prop}\label{ld0coeffs}
	In the setting of proposition \ref{dicoeffs}, the generating classes $\lambda,\,\delta_{0}^{\rm t},\,\delta_{0}^{\rm p},\,\delta_{0}^{\rm b}\in\Pic(\OK{R}_{g})_{\B{Q}}$ have coefficients
	\[
	\begin{array}{rc|crclcl}
		\multirow{4}{*}[-7pt]{$\varrho^{+}_{\rm null}$}
		&&\;\!&\lambda^{+}
		&=&2^{g-3}(2^{g-1}+1)&=&2^{g-3}(2^{g-1}+1)
		\\[2mm]
		&&&\delta_{0}^{\rm t, +}
		&=&0&=&0
		\\[2mm]
		&&&\delta_{0}^{\rm p, +}
		&=&2^{2g-7}&=&2^{g-3}\,2^{-2}\,2^{g-2}
		\\[2mm]
		&&&\delta_{0}^{\rm b, +}
		&=&2^{g-5}(2^{g-1}+1)&=&2^{g-3}\,2^{-2}(2^{g-1}+1)
		\\[-1mm]
		\multicolumn{1}{c}{}
		\\
		\multirow{4}{*}[-7pt]{$\varrho^{-}_{\rm null}$}
		&&&\lambda^{-}
		&=&2^{2g-4}&=&2^{g-3}\,2^{g-1}
		\\[2mm]
		&&&\delta_{0}^{\rm t, -}
		&=&2^{2g-6}&=&2^{g-3}\,2^{-2}\,2^{g-1}
		\\[2mm]
		&&&\delta_{0}^{\rm p, -}
		&=&2^{2g-7}&=&2^{g-3}\,2^{-2}\,2^{g-2}
		\\[2mm]
		&&&\delta_{0}^{\rm b, -}
		&=&2^{g-5}(2^{g-1}-1)&=&2^{g-3}\,2^{-2}(2^{g-1}-1)
	\end{array}
	\]
	in the rational expansions of the Prym-null classes in genus $g$.
\end{prop}
\begin{dem}
	Since the $\delta_{1}^{\rm n}$ coefficients have already been computed (proposition \ref{dicoeffs}), it is straightforward to check that the linear relation provided by the family $Y_{0}$ in each case directly determines the corresponding $\delta_{0}^{\rm t}$ coefficient:
	\[
	\begin{array}{lcl}
		Y_{0}&\rightsquigarrow&
		\left\lbrace
		\begin{array}{lcl}
			\delta_{0}^{\rm t, +}&=&0
			\vspace*{2mm}\\
			\delta_{0}^{\rm t, -}&=&2^{2g-6}
		\end{array}\right.
	\end{array}
	\]
	Plugging these into lemma \ref{LIlinrel}, we obtain the following linear systems:
	\[
	\begin{array}{c}
		\left(\begin{array}{ccc}
			1&-12&\phantom{-}0\\
			1&\phantom{-}0&-4\\
			1&-4&-4
		\end{array}\right)
		\cdot
		\left(\begin{array}{c}
			\lambda^{+}\\
			\delta_{0}^{\rm p, +}\\
			\delta_{0}^{\rm b, +}
		\end{array}\right)
		=
		\left(\begin{array}{l}
			-2^{g-3}(2^{g-2}-1)\\
			\phantom{-}0\\
			-2^{2g-5}
		\end{array}\right)
		\vspace*{5mm}\\
		\left(\begin{array}{lll}
			1&-12&\phantom{-}0\\
			1&\phantom{-}0&-4\\
			1&-4&-4
		\end{array}\right)
		\cdot
		\left(\begin{array}{c}
			\lambda^{-}\\
			\delta_{0}^{\rm p, -}\\
			\delta_{0}^{\rm b, -}
		\end{array}\right)
		=
		\left(\begin{array}{l}
			-2^{2g-5}\\
			\phantom{-}2^{g-3}\\
			-2^{g-3}(2^{g-2}-1)
		\end{array}\right)
	\end{array}
	\]
	The two sets of solutions are precisely the expressions stated above. \qed
\end{dem}

%% file: 235ClassDnull.tex
For the first time, all of the rational coefficients introduced in remark \ref{defPnullcoeffs} are known to us, by virtue of propositions \ref{dicoeffs} and \ref{ld0coeffs}. As a result, we are finally in a position to express the rational classes of $\OK{P}_{\rm\! null}^{+}$ and $\OK{P}_{\rm\! null}^{-}$ in terms of the generating classes of $\Pic(\OK{R}_{g})_{\B{Q}}$, which was our main goal.

\begin{thm}\label{classprymnulls}
	For $g\geq 5$, the classes $\varrho^{+}_{\rm null},\,\varrho^{-}_{\rm null}\in\Pic(\OK{R}_{g})_{\B{Q}}$ are given by
	\[
	\begin{array}{rccclc}
		\varrho^{+}_{\rm null}
		&=&\multicolumn{4}{l}{
			2^{g-3}\,\bigg(
			(2^{g-1}+1)\,\lambda
			\hspace*{5.0pt}-\hspace*{5.0pt}
			\dfrac{1}{4}\,\Big(
			2^{g-2}\,\delta_{0}^{\rm p}+
			(2^{g-1}+1)\,\delta_{0}^{\rm b}
			\Big)}
		\vspace*{1.5mm}\\
		&&&-&
		\displaystyle\sum_{i=1}^{k}\,
		\Big(
		(2^{i-1}-1)(2^{g-i}-1)\,\delta_{i}^{\rm n}+
		(2^{i}-1)(2^{g-i-1}-1)\,\delta_{i}^{\rm t}
		\hspace*{3.5pt}+
		\vspace*{1mm}\\
		&&&&\multicolumn{1}{r}{+\hspace*{3.5pt}
			(2^{g-1}-2^{i-1}-2^{g-i-1}+1)\,\delta_{i}^{\rm p}
			\Big)}
		\vspace*{3.5mm}\\
		&&&-&\multicolumn{2}{r}{\psi(g)\cdot\Big(
			(2^{g/2-1}-1)(2^{g/2}-1)\,\delta_{g/2}^{\rm n}+
			(2^{g-1}-2^{g/2}+1)\,\delta_{g/2}^{\rm p}
			\Big)\bigg)}
		\vspace*{5mm}\\
		\varrho^{-}_{\rm null}
		&=&\multicolumn{4}{l}{
			2^{g-3}\,\bigg(
			2^{g-1}\,\lambda
			\hspace*{5.0pt}-\hspace*{5.0pt}
			\dfrac{1}{4}\,\Big(
			2^{g-1}\,\delta_{0}^{\rm t}+
			2^{g-2}\,\delta_{0}^{\rm p}+
			(2^{g-1}-1)\,\delta_{0}^{\rm b}
			\Big)}
		\vspace*{1.5mm}\\
		&&&-&
		\displaystyle\sum_{i=1}^{k}\,
		\Big(
		2^{i-1}\,(2^{g-i}-1)\,\delta_{i}^{\rm n}+
		(2^{i}-1)\,2^{g-i-1}\,\delta_{i}^{\rm t}
		\hspace*{3.5pt}+
		\vspace*{1mm}\\
		&&&&\multicolumn{1}{r}{+\hspace*{3.5pt}
			(2^{g-1}-2^{i-1}-2^{g-i-1})\,\delta_{i}^{\rm p}
			\Big)}
		\vspace*{3.5mm}\\
		&&&-&\multicolumn{2}{l}{\psi(g)\cdot\Big(
			2^{g/2-1}\,(2^{g/2}-1)\,\delta_{g/2}^{\rm n}+
			(2^{g-1}-2^{g/2})\,\delta_{g/2}^{\rm p}
			\Big)\bigg)}
	\end{array}
	\]
	where the upper bound $k$ and the parity-checking function $\psi(g)$, defined as
	\[
	\begin{array}{rclcl}
		k&=&\left\lceil g/2 \right\rceil-1&=&
		\left\lbrace
		\begin{array}{ll}
			\left\lfloor g/2 \right\rfloor&\textrm{\rm{ if $g$ odd}}\\[1mm]
			\left\lfloor g/2 \right\rfloor-1&\textrm{\rm{ if $g$ even}}
		\end{array}\right.
		\\[6mm]
		\psi(g)&=&\dfrac{1+(-1)^{g}}{2}&=&
		\left\lbrace
		\begin{array}{ll}
			\mathrlap{0}\phantom{\left\lfloor g/2 \right\rfloor-1}&\textrm{\rm{ if $g$ odd}}	\\[1mm]
			1&\textrm{\rm{ if $g$ even}}
		\end{array}\right.
	\end{array}
	\]
	account for the slight variation that occurs when $g=2i$.
\end{thm}

\begin{exam}[quartic tails]\label{testquartail}
	We fix a general curve $(C,p)\in\K{M}_{g-3,\,1}$ and a general pencil $\gamma\colon\N{Bl}_{16}(\B{P}^2)\to\B{P}^{1}$ of plane quartics, with fibers
	\[
	\begin{array}{rcl}
		\{Q_{\lambda}=\gamma^{-1}(\lambda)\}_{\lambda\in\B{P}^1}&\subset&\OK{M}_{3}
	\end{array}
	\]
	together with a section $\zeta\colon\B{P}^{1}\to\N{Bl}_{16}(\B{P}^2)$ induced by one of the basepoints. We may then glue the curve $(C,p)$ to the pencil $\gamma$ along $\zeta$, thus producing a pencil of stable curves
	\[
	\K{Q}=(C\times \B{P}_{1})\cup_{\{p\}\times\B{P}^{1}\:\!\sim\, \zeta(\B{P}^{1})}\N{Bl}_{16}(\B{P}^2)\lto \B{P}^{1}
	\]
	which corresponds to
	\[
	\begin{array}{rclclcl}
		\K{Q}&\equiv&\{C\cup_{p\sim\zeta(\lambda)}Q_{\lambda}\}_{\lambda\in\B{P}^{1}}&\subset&\Delta_{3}&\subset&\OK{M}_{g}
	\end{array}
	\]
	Standard techniques show its intersection table to be:
	\[
	\begin{array}{c|cccccccc}
		&\lambda&\delta_{0}&\delta_{1}&\delta_{2}&\delta_{3}&\delta_{(j\geq 4)}\\[\dimexpr-\normalbaselineskip+0.5mm]
		\\\hline
		\\[\dimexpr-\normalbaselineskip+1mm]
		\K{Q}&3&27&0&0&-1&0
	\end{array}
	\]
	We now fix a Prym root $\eta_C\in R_{g-3}(C)$ and lift $\K{Q}$ to a test curve $R$, as follows:
	\[
	\begin{array}{rclclcl}
		R&\equiv&\{(C\cup_{p\sim\zeta(\lambda)}Q_{\lambda},\,(\eta_C,\K{O}_{Q_{\lambda}}))\}_{\lambda\in\B{P}^1}&\subset&\Delta_{3}^{\rm t}&\subset&\OK{R}_g
	\end{array}
	\]
	Observe that $\pi_{*}(R)=\K{Q}$. Then $R\cdot\lambda=\K{Q}\cdot\lambda=3$ and
	\[\begin{array}{rclcl}
		R\cdot\delta_{3}^{\rm t}
		&=&\K{Q}\cdot\delta_{3}&=&-1
	\end{array}
	\]
	If we look at the $27$ points $\lambda_{\infty}\in\B{P}^1$ corresponding to singular quartics of $\gamma$ and blow up the node of the component $Q_{\lambda_{\infty}}\in\Delta_{0}$, we can see that the pullback of $\eta_{\lambda_{\infty}}=(\eta_C,\K{O}_{Q_{\lambda_{\infty}}})$ is $(\eta_C,\K{O}_{\B{P}^1})$, which is nontrivial. Hence $R_{\lambda_{\infty}}\in\Delta_{0}^{\rm p}$ and
	\[
	\begin{array}{rclcl}
		R\cdot \delta_{0}^{\rm p}
		&=&\K{Q}\cdot\delta_{0}&=&27
	\end{array}
	\]
	All other intersection numbers are $0$, so we get a table:
	\[
	\begin{array}{c|cccccccc}
		&\lambda&\delta_{0}^{\rm t}&\delta_{0}^{\rm p}&\delta_{0}^{\rm b}&\delta_{3}^{\rm n}&\delta_{3}^{\rm t}&\delta_{3}^{\rm p}&\delta_{(j\neq 0,\,3)}\\[\dimexpr-\normalbaselineskip+0.5mm]
		\\\hline
		\\[\dimexpr-\normalbaselineskip+1mm]
		R&3&0&27&0&0&-1&0&0
	\end{array}
	\]
	Note the similarities between $R$ and the family $F_{0}$ from example \ref{testelliptail}.
\end{exam}

\begin{rem}\label{interquartail}
	Applying theorem \ref{classprymnulls} to example \ref{testquartail}, we get
	\[
	\begin{array}{rcccl}
		R\cdot\OK{P}_{\rm\! null}^{+}&=&3\,\lambda^{+}-27\,\delta_{0}^{\rm p, +}+\delta_{3}^{\rm t, +}&=&2^{g-1}(2^{g-4}-1)
		\\[2mm]
		R\cdot\OK{P}_{\rm\! null}^{-}&=&3\,\lambda^{-}-27\,\delta_{0}^{\rm p, -}+\delta_{3}^{\rm t, -}&=&2^{2g-5}
	\end{array}
	\]
	These intersection numbers may in fact be interpreted, as the limit linear series techniques introduced in earlier cases are also quite useful here.
\end{rem}

Again, we want to consider $R\cap\OK{P}_{\rm\! null}^{+}$ (resp. $\OK{P}_{\rm\! null}^{-}$). If a stable Prym curve
\[
\begin{array}{rclcl}
	R_{\lambda}&=&(C\cup_{p\sim z}Q_{\lambda},\,(\eta_C,\K{O}_{Q_{\lambda}}))&\in&R
\end{array}
\]
lies in $\OK{P}_{\rm\! null}^{+}$ (resp. $\OK{P}_{\rm\! null}^{-}$), we can produce a limit $\mathfrak{g}_{g-1}^{1}$ on $C\cup_{p\sim z}Q_{\lambda}$ such that
\[
\left\lbrace
\begin{array}{lcl}
	L_C&=&\theta_C\otimes\eta_C\,(3p)
	\\[2mm]
	L_{Q_{\lambda}}&=&\theta_{Q_{\lambda}}\,((g-3)z)
\end{array}
\right.
\]
with $h^{0}(\theta_C\otimes\eta_C)+h^{0}(\theta_{Q_{\lambda}})\equiv 0 \mod 2$, where $\theta_C$ and $\theta_{Q_{\lambda}}$ have the same parity (resp. opposite parity). Then $\theta_C\otimes\eta_C$ and $\theta_{Q_{\lambda}}$ have the same parity and we get the following possibilities:
\[
\begin{array}{cccccccl}
	&
	&h^{0}(\theta_C)&h^{0}(\theta_{Q_{\lambda}})&h^{0}(\theta_C\otimes\eta_C)
	\\[\dimexpr-\normalbaselineskip+0.5mm]
	\\\cline{2-5}
	&\\[\dimexpr-\normalbaselineskip+1mm]
	\multirow{2}{*}{$\OK{P}_{\rm\! null}^{+}$}
	&&0&0&0&&\rightsquigarrow&\textrm{contradiction}\\
	&&1&1&1&&\rightsquigarrow&\N{(R,+,1)}
	\\[\dimexpr-\normalbaselineskip+0mm]
	&\\\cline{2-5}
	&\\[\dimexpr-\normalbaselineskip+1mm]
	\multirow{2}{*}{$\OK{P}_{\rm\! null}^{-}$}
	&&1&0&0&&\rightsquigarrow&\textrm{contradiction}\\
	&&0&1&1&&\rightsquigarrow&\N{(R,-,1)}
	\\[\dimexpr-\normalbaselineskip+0mm]
	&\\\cline{2-5}
\end{array}
\]

In order to deal with the cases not covered by remark \ref{contradict}, we first need to understand how theta characteristics on canonical genus $3$ curves look like. Let us quickly elaborate on this.

Given a nonhyperelliptic arbitrary curve $X\in\K{M}_3$, its canonical embedding realises it as a plane quartic $Q\inj\B{P}^2$, with the canonical series manifesting as the restriction of the hyperplane series to the curve.

Take a theta characteristic $\theta$ on $X$, with $\theta^{\otimes 2}\cong\omega_{X}\in W_{4}^{2}(X)$. Then we have $\deg(\theta)=2$, and $\theta$ is of type $\G{g}_{2}^{r}$ on $X$ whenever $h^{0}(\theta)=r+1>0$. But $X$ is not hyperelliptic, so it does not admit any $\G{g}_{2}^{1}$ and thus $r\leq 0\,\Rightarrow\, h^{0}(\theta)=r+1\leq 1$. Therefore, the $36$ even theta characteristics of $X$ have $h^{0}(\theta)=0$ and the $28$ odd ones have $h^{0}(\theta)=1$. In particular, $X$ has no vanishing theta-nulls.

Let $\theta$ be odd. Then $\vert\theta\vert=\{D\}$ with $2D\sim K_X$, that is, $D=x+y$ and
\[
\begin{array}{rclcl}
	2D&=&2x+2y&=&H\cap Q
\end{array}
\]
for some hyperplane $H\inj\B{P}^2$. If moreover $x\neq y$ (for example, for $X$ general), we get a one-to-one correspondence between odd theta characteristics of $X$ and bitangents to its canonical model $Q$. Note that, if $X$ is special enough for $Q$ to have any hyperflexes, then the tangent lines at such points must be included in the correspondence too.

\begin{posss}[$\N{R,+,1)}$, $\N{(R,-,1}$]
	In both of these cases, we have
	\[
	\begin{array}{rcl}
		h^{0}(\theta_C\otimes\eta_C)=1&\Rightarrow&a_{0}^{\ell_C}(p)<3\leq a_{1}^{\ell_C}(p)
		\vspace*{2mm}\\
		h^{0}(\theta_{Q_{\lambda}})=1&\Rightarrow&a_{0}^{\ell_{\smash{Q_{\lambda}}}}(z)<g-3\leq a_{1}^{\ell_{\smash{Q_{\lambda}}}}(z)
	\end{array}
	\]
	Now $(C,p)$ is general, so we may assume that $p\notin\N{supp}(\theta_C\otimes\eta_C)$. Therefore
	\[\begin{array}{rcl}
		h^{0}(\theta_C\otimes\eta_C\,(-p))=h^{0}(\theta_C\otimes\eta_C)-1=0&\Rightarrow&a_{1}^{\ell_C}(p)<4
		\vspace*{2mm}\\
		h^{0}(\theta_C\otimes\eta_C\,(p))=h^{0}(\theta_C\otimes\eta_C)=1&\Rightarrow&a_{0}^{\ell_C}(p)<2\leq a_{1}^{\ell_C}(p)
	\end{array}\]
	Hence $a_{1}^{\ell_C}(p)=3$ and, via the limit $\G{g}_{g-1}^{1}$ condition, $a_{0}^{\ell_{\smash{Q_{\lambda}}}}(z)=g-4$. Moreover, we may assume that the basepoint $z$ is not a hyperflex of $Q_{\lambda}$, as the pencil $\gamma$ is also general. Consequently, $\N{supp}(\theta_{Q_{\lambda}})$ does not consist of $z$ twice, that is,
	\[\begin{array}{rcl}
		\N{div}(\theta_{Q_{\lambda}})\neq 2z&\Rightarrow&a_{1}^{\ell_{\smash{Q_{\lambda}}}}(z)\leq g-2
	\end{array}\]
	which combined with the condition $a_{0}^{\ell_C}(p)+a_{1}^{\ell_{\smash{Q_{\lambda}}}}(z)\geq g-1$ yields $a_{0}^{\ell_C}(p)=1$ and $a_{1}^{\ell_{\smash{Q_{\lambda}}}}(z)=g-2$. In turn, this means that $z\in\N{supp}(\theta_{Q_{\lambda}})$, and that $\ell$ is a refined limit $\G{g}_{g-1}^{1}$ of the form
	\[\begin{array}{cclcl}
		\ell_C&=&\vert\theta_C\otimes\eta_C\,(2p)\vert+p&\in&G_{g-1}^{1}(C)
		\vspace*{2mm}\\
		\ell_{Q_{\lambda}}&=&\vert\theta_{Q_{\lambda}}\,(z)\vert+(g-4)z&\in&G_{g-1}^{1}(Q_{\lambda})
	\end{array}\]
	with vanishing sequences $(1,\,3)$ and $(g-4,\,g-2)$.
\end{posss}

In conclusion, for each pair $(Q_{\lambda},\,\theta_{Q_{\lambda}})$ consisting of a plane quartic $Q_{\lambda}$ of $\gamma$ equipped with an odd theta characteristic $\theta_{Q_{\lambda}}$ such that $z=\zeta(\lambda)\in\N{supp}(\theta_{Q_{\lambda}})$, then every $\theta_C\in S_{g-3}^{-}(C)$ with $\theta_C\otimes\eta_C\in S_{g-3}^{-}(C)$ yields a limit $\G{g}_{g-1}^{1}$ as above, and these limit linear series are the only ones contributing to the intersection $R\cap\OK{P}_{\rm\! null}^{+}$ (resp. $\theta_C\in S_{g-3}^{+}(C)$ with $\theta_C\otimes\eta_C\in S_{g-3}^{-}(C)$, $R\cap\OK{P}_{\rm\! null}^{-}$). The natural question then arises as to how many such pairs $(Q_{\lambda},\,\theta_{Q_{\lambda}})$ there are.

We have discussed that odd theta characteristics $\theta_{Q_{\lambda}}$ of the plane curve $Q_{\lambda}$ correspond to bitangents to the quartic. Under this identification, the condition $z=\zeta(\lambda)\in\N{supp}(\theta_{Q_{\lambda}})$ corresponds to the bitangent having the basepoint $z$ as one of its contact points. In particular, for each $Q_{\lambda}$ we get only one candidate: the tangent line $T_{z}(Q_{\lambda})\subset\B{P}^2$, which will intersect $Q_{\lambda}$ in two additional points. If we find out for how many values of $\lambda$ these two points coincide, we will have found the pairs $(Q_{\lambda},\,\theta_{Q_{\lambda}})\equiv(Q_{\lambda},T_{z}(Q_{\lambda}))$ we are trying to count.

Now, we can study the pencil $\gamma$ by taking two general polynomials
\[
\begin{array}{cccc}
	F(\F{x})=\sum_{i+j+k=4}\,a_{ijk}\,x_0^{i}x_1^{j}x_2^{k}\;,& G(\F{x})=\sum_{i+j+k=4}\,b_{ijk}\,x_0^{i}x_1^{j}x_2^{k}&\in&\B{C}[x_0,x_1,x_2]_{4}
\end{array}
\]
and considering the family $\{Q_{\lambda}\}$ described by $H(\F{x},\lambda)=\lambda_0\,F(\F{x})+\lambda_1\,G(\F{x})=0$, with basepoints $\{(\F{x},\lambda)\;\slash\;F(\F{x})=G(\F{x})=0\}\ni\zeta(\lambda)=z$. By a suitable change of coordinates, we may assume $z=(1:0:0)\in\B{P}^2$. If we write $H(\F{x},\lambda)$ as
\[
\begin{array}{ccccccc}
	H(\F{x},\lambda)&=&H_{\lambda}(\F{x})&=&\sum_{i+j+k=4}\,c_{ijk}(\lambda)\,x_0^{i}x_1^{j}x_2^{k}&\in&\B{C}[x_0,x_1,x_2]_{4}
\end{array}
\]
with $c_{ijk}(\lambda)=\lambda_0\,a_{ijk}+\lambda_1\,b_{ijk}\in\B{C}[\lambda]_{1}$, this means that $c_{400}(\lambda)=0$.

Moreover, the tangent line $T_{z}(Q_{\lambda})$ is given by
\[
\begin{array}{ccccc}
	\frac{\partial H_{\lambda}(\F{x})}{\partial x_0}(z)\,x_0+\frac{\partial H_{\lambda}(\F{x})}{\partial x_1}(z)\,x_1+\frac{\partial H_{\lambda}(\F{x})}{\partial x_2}(z)\,x_2&=&c_{310}(\lambda)\,x_1+c_{301}(\lambda)\,x_2&=&0
\end{array}
\]
Switching to coordinates $u,v$ on $T_{z}(Q_{\lambda})=\B{P}^{1}_{u,v}$, that is,
\[
\begin{array}{ccc}
	x_0=u\;,&x_1=-c_{301}(\lambda)\,v\;,&x_2=c_{310}(\lambda)\,v
\end{array}
\]
we see that $z=\{v=0\}=(1:0)\in\B{P}^1_{u,v}$ and that the intersection $Q_{\lambda}\cap T_{z}(Q_{\lambda})$ is given by
\[
\begin{array}{ccccc}
	\OC{H}_{\lambda}(u,v)&=&\sum_{i+j+k=4}\,(-1)^{j}\,c_{301}(\lambda)^{j}\,c_{310}(\lambda)^{k}\,c_{ijk}(\lambda)\,u^{i}v^{j+k}&=&0
\end{array}
\]
Since the intersection contains $z$ twice, this polynomial has no $v^{0}$, $v^{1}$ terms:
\[
\begin{array}{ccccc}
	j+k=0&\Rightarrow&i=4&\rightsquigarrow&c_{400}(\lambda)=0
	\vspace*{2mm}\\
	j+k=1&\Rightarrow&i=3&\rightsquigarrow&-c_{301}(\lambda)\,c_{310}(\lambda)+c_{310}(\lambda)\,c_{301}(\lambda)=0
\end{array}
\]
Factoring out $v^2$, we get a quadric
\[
\begin{array}{ccl}
	q_{\lambda}(u,v)&=&
	\sum_{
		\begin{subarray}{l}
			0\leq i\leq 2\vspace*{.25mm}\\
			j+k=4-i
	\end{subarray}}\,(-1)^{j}\,c_{301}(\lambda)^{j}\,c_{310}(\lambda)^{k}\,c_{ijk}(\lambda)\,u^{i}v^{2-i}\vspace*{2mm}\\
	&=&c^{uu}_{3}(\lambda)\,u^{2}+c^{uv}_{4}(\lambda)\,uv+c^{vv}_{5}(\lambda)\,v^{2}
\end{array}
\]
whose roots correspond to the two additional points lying in $Q_{\lambda}\cap T_{z}(Q_{\lambda})$. The degree of each summand $(-1)^{j}\,c_{301}(\lambda)^{j}\,c_{310}(\lambda)^{k}\,c_{ijk}(\lambda)$ is $j+k+1$, so we have
\[
\begin{array}{ccccc}
	c^{uu}_{3}(\lambda)&=&\sum_{j+k=2}\,(-1)^{j}\,c_{301}(\lambda)^{j}\,c_{310}(\lambda)^{k}\,c_{2jk}(\lambda)&\in&\B{C}[\lambda_0,\lambda_{1}]_{3}\vspace*{2mm}\\
	c^{uv}_{4}(\lambda)&=&\sum_{j+k=3}\,(-1)^{j}\,c_{301}(\lambda)^{j}\,c_{310}(\lambda)^{k}\,c_{1jk}(\lambda)&\in&\B{C}[\lambda_0,\lambda_{1}]_{4}\vspace*{2mm}\\
	c^{vv}_{5}(\lambda)&=&\sum_{j+k=4}\,(-1)^{j}\,c_{301}(\lambda)^{j}\,c_{310}(\lambda)^{k}\,c_{0jk}(\lambda)&\in&\B{C}[\lambda_0,\lambda_{1}]_{5}
\end{array}
\]
Finally, the values of $\lambda$ for which the roots of $q_{\lambda}(u,v)$ coincide are determined by the roots of the discriminant
\[
\begin{array}{ccccccc}
	\Delta(\lambda)&=&\Delta(q_{\lambda}(u,v))&=&c^{uv}_{4}(\lambda)^{2}-4\,c^{uu}_{3}(\lambda)\,c^{vv}_{5}(\lambda)&\in&\B{C}[\lambda_0,\lambda_{1}]_{8}
\end{array}
\]
which is an octic polynomial. Therefore we obtain
\[
\begin{array}{ccccccc}
	\#\{(Q_{\lambda},\theta_{Q_{\lambda}})\;\slash\;z\in\N{supp}(\theta_{Q_{\lambda}})\}&=&\#\{\lambda\in\B{P}^{1}\;\slash\;\Delta(\lambda)=0\}&=&8&=&2^3
\end{array}
\]
and we can finally observe the appearance of the intersection numbers provided by theorem \ref{classprymnulls} and remark \ref{interquartail}, as the count becomes:
\[
\begin{array}{lclcl}
	\#\{(Q_{\lambda},\theta_{Q_{\lambda}})\;\slash\;z\in\N{supp}(\theta_{Q_{\lambda}})\}\;\cdot
	\\[1mm]
	\cdot\;\#\{\theta_C\in S_{g-3}^{-}(C)\;\slash\;\theta_C\otimes\eta_C\in S_{g-3}^{-}(C)\}&=&2^{3}\cdot\N{N}_{g-3}^{-}&=&2^{g-1}(2^{g-4}-1)
	\vspace*{4mm}\\
	\#\{(Q_{\lambda},\theta_{Q_{\lambda}})\;\slash\;z\in\N{supp}(\theta_{Q_{\lambda}})\}\;\cdot
	\\[1mm]
	\cdot\;\#\{\theta_C\in S_{g-3}^{+}(C)\;\slash\;\theta_C\otimes\eta_C\in S_{g-3}^{-}(C)\}&=&2^{3}\cdot\N{N}_{g-3}^{\pm}&=&2^{2g-5}
\end{array}
\]
In particular, this indicates the lack of contribution from singular fibers.

\begin{exam}[more irreducible nodal curves]\label{moretestirrednod}
	If we recall the family
	\[
	\begin{array}{rclclcl}
		\K{Y}&\equiv&\{B_{py}\}_{y\in B}&\subset&\Delta_{0}&\subset&\OK{M}_{g}
	\end{array}
	\]
	from example \ref{testirrednod}, which was lifted to a test curve $Y_{0}\subset\Delta_{0}^{\rm t}\subset\OK{R}_g$, then there are two more standard lifts $Z_{0}$ and $T_{0}$ in $\OK{R}_{g}$, which arise when $\K{Y}$ is pulled back by the maps $\Delta_{0}^{\rm p}\to\Delta_{0}$ and $\Delta_{0}^{\rm b}\to\Delta_{0}$ respectively:
	\[
	\begin{array}{rclclcl}
		Z_0&\equiv&\{(B_{py},\,\eta_{y}^{\rm p})\;\slash\;\eta_{y}^{\rm p}\in\Delta_{0}^{\rm p}(B_{py}) \}_{y\in B}&\subset&\Delta_{0}^{\rm p}&\subset&\OK{R}_g
		\\[2mm]
		T_0&\equiv&\{(B\cup_{p\sim 0,\,y\sim\infty}E,\,\eta_{y}^{\rm b})\;\slash\;\eta_{y}^{\rm b}\in\Delta_{0}^{\rm b}(B_{py}) \}_{y\in B}&\subset&\Delta_{0}^{\rm b}&\subset&\OK{R}_g
	\end{array}
	\]
	If we set $k=\#R_{g-1}(B)=2^{2g-2}-1$, we can see that their intersection table is:
	\[
	\begin{array}{c|cccccccc}
		&\lambda&\delta_{0}^{\rm t}&\delta_{0}^{\rm p}&\delta_{0}^{\rm b}&\delta_{1}^{\rm n}&\delta_{1}^{\rm t}&\delta_{1}^{\rm p}&\delta_{(j\geq 2)}\\[\dimexpr-\normalbaselineskip+0.5mm]
		\\\hline
		\\[\dimexpr-\normalbaselineskip+1mm]
		Y_0&0&2-2g&0&0&1&0&0&0\\
		Z_0&0&0&4k(1-g)&0&0&k&k&0\\
		T_0&0&0&0&2^{2g-2}(1-g)&1&0&k&0
	\end{array}
	\]
	Note that $\deg(\Delta_{0}^{\rm p}\vert\Delta_{0})=2k$ and $\deg(\Delta_{0}^{\rm b}\vert\Delta_{0})=2^{2g-2}=k+1$.
\end{exam}

\begin{rem}\label{interirrednod}
	Applying theorem \ref{classprymnulls} to example \ref{moretestirrednod}, it follows that
	\[
	\begin{array}{rclcl}
		Z_{0}\cdot\OK{P}_{\rm\! null}^{+}&=&(2^{2g-2}-1)\,\mu
		&=&\#R_{g-1}(B)\cdot\mu
		\\[2mm]
		Z_{0}\cdot\OK{P}_{\rm\! null}^{-}&=&(2^{2g-2}-1)\,\mu
		&=&\#R_{g-1}(B)\cdot\mu
		\\[4mm]
		T_{0}\cdot\OK{P}_{\rm\! null}^{+}&=&2^{g-2}(2^{g-1}+1)\,\mu
		&=&\#S_{g-1}^{+}(B)\cdot\mu
		\\[2mm]
		T_{0}\cdot\OK{P}_{\rm\! null}^{-}&=&2^{g-2}(2^{g-1}-1)\,\mu
		&=&\#S_{g-1}^{-}(B)\cdot\mu
	\end{array}
	\]
	with the factor $\mu=Y_0^{\rm n}\cdot\OC{\Theta}_{\rm null}=2^{g-3}(2^{g-2}(g-3)+1)$ indicating the number of nodal curves $B_{py}$ in $\K{Y}$ that admit a vanishing theta-null $\theta_{y}\in\OC{\Theta}_{\rm null}(B_{py})$, which we computed in the argument preceding preposition \ref{ld0coeffs}. Once more, it may be interesting to provide an interpretation of these results.
\end{rem}

According to example \ref{BoundDivRgO}, any Prym root $\eta_{B}\in R_{g-1}(B)$ gives rise to two elements $\eta_{y}^{\rm p, +},\,\eta_{y}^{\rm p, -}\in\Delta_{0}^{\rm p}(B_{py})$, depending on which of the two possible gluings $\restr{\eta_{B}}{p}\cong\restr{\eta_{B}}{y}$ is chosen. In particular, for each pair $(B_{py},\,\theta_{y})\in\OC{\Theta}_{\rm null}$, tensoring $\theta_{y}$ by either $\eta_{y}^{\rm p, +}$ or $\eta_{y}^{\rm p, -}$ produces stable spin curves of opposite parity, so that
\[
\begin{array}{rcl}
	(B_{py},\,\eta_{y}^{\rm p, +})&\in& Z_{0}\cap\OK{P}_{\rm\! null}^{+}
	\\[2mm]
	(B_{py},\,\eta_{y}^{\rm p, -})&\in& Z_{0}\cap\OK{P}_{\rm\! null}^{-}
\end{array}
\]
which explains the emergence of the factors
\[
\begin{array}{rcl}
	k&=&\#\{\eta_{B}\in R_{g-1}(B)\}
	\\[2mm]
	\mu&=&\#\{(B_{py},\,\theta_{y})\in\OC{\Theta}_{\rm null}\}
\end{array}
\]
in the intersection numbers $Z_{0}\cdot\OK{P}_{\rm\! null}^{+}$ and $Z_{0}\cdot\OK{P}_{\rm\! null}^{-}$.

Similarly, for each pair $(B_{py},\,\theta_{y})\in\OC{\Theta}_{\rm null}$, the root $\restr{\theta_{y}}{B}\in\sqrt{\omega_{B}(p+q)}$ can be subtracted from any theta characteristic $\theta_{B}\in S_{g-1}(B)=\sqrt{\omega_{B}}\hspace*{1pt}$ so as to create a root $\eta_{B}\in\sqrt{\K{O}_{B}(-p-q)}$. This in turn yields a unique stable Prym curve
\[
\begin{array}{rclcl}
	(X,\eta_{y}^{\rm b})&=&(B\cup_{p\sim 0,\,y\sim\infty}E,\,\eta_{y}^{\rm b})&\in&T_{0}\cap\OK{P}_{\rm\! null}
\end{array}
\]
such that $\eta_{y}^{\rm b}$ restricts to $(\eta_{B},\,\K{O}_{E}(1))$ on $\Pic(B)\oplus\Pic(E)$. Furthermore, $(X,\eta_{y}^{\rm b})$ lies in $\OK{P}_{\rm\! null}^{+}$ (resp. $\OK{P}_{\rm\! null}^{-}$) whenever $\theta_{B}$ is even (resp. odd) by construction of the Prym-null divisors, which brings to light the connection between
\[
\begin{array}{rclcl}
	\#S_{g-1}^{+}(B),&\#S_{g-1}^{-}(B),&\mu
\end{array}
\]
and the intersection numbers $T_{0}\cdot\OK{P}_{\rm\! null}^{+}$ and $T_{0}\cdot\OK{P}_{\rm\! null}^{-}$.